\documentclass[reqno,12pt,letterpaper]{article}
\pdfoutput=1
\usepackage{color, xcolor, colortbl}
\usepackage{graphicx,epstopdf}
\usepackage{geometry}
\usepackage{enumerate}
\usepackage{amsmath,amssymb,amsthm}
\usepackage{algorithm}
\usepackage{algorithmic}
\usepackage{caption}
\usepackage{comment}
\usepackage{subcaption}
\usepackage{bm}
\usepackage{appendix}
\usepackage{multirow}
\usepackage{braket}
\usepackage[english]{babel}
\usepackage{hyperref}
\usepackage[capitalize]{cleveref}
\usepackage[sort&compress,square,numbers]{natbib}
\usepackage{adjustbox}
\usepackage{xspace}
\usepackage[roman]{complexity}
\usepackage{soul}
\usepackage{tikz}
\usetikzlibrary{quantikz}
\usepackage{rotating}
\usepackage{setspace}
\usepackage{fancyhdr}

\newcommand{\dd}{\partial}
\newcommand{\la}{\langle}
\newcommand{\ra}{\rangle}

\usepackage{authblk}

\newcommand{\norm}[1]{\left\lVert#1\right\rVert}

\newcommand{\op}{{\rm op}}

\newcommand{\Or}{\mathcal{O}}

\newcommand{\RR}{\mathbb{R}}

\newcommand{\TT}{\mathbb{T}}

\DeclareMathOperator*{\cinf}{C^{\infty}}

\newcommand{\ad}{\operatorname{ad}}

\newtheorem{thm}{\protect\theoremname}
\theoremstyle{plain}
\newtheorem{lemma}[thm]{\protect\lemmaname}
\theoremstyle{plain}

\theoremstyle{plain}
\theoremstyle{plain}

\theoremstyle{plain}

\newcommand{\Op}{\operatorname{op}}

\providecommand{\definitionname}{Definition}
\providecommand{\assumptionname}{Assumption}
\providecommand{\corollaryname}{Corollary}
\providecommand{\lemmaname}{Lemma}
\providecommand{\propositionname}{Proposition}
\providecommand{\remarkname}{Remark}
\providecommand{\theoremname}{Theorem}

\newcommand\blfootnote[1]{%
  \begingroup
  \renewcommand\thefootnote{}\footnote{#1}%
  \addtocounter{footnote}{-1}%
  \endgroup
}

\newcommand{\REV}[1]{#1}

\title{Discrete Superconvergence Analysis for Quantum Magnus Algorithms of Unbounded Hamiltonian Simulation}

\author[1]{Yonah Borns-Weil}
\author[2,3]{Di Fang}
\author[2]{Jiaqi Zhang}

\affil[1]{Department of Mathematics, University of California, Berkeley}
\affil[2]{Department of Mathematics, Duke University}
\affil[3]{Duke Quantum Center, Duke University}

\date{} 

\begin{document}
\maketitle

\begin{abstract}
Motivated by various applications, unbounded Hamiltonian simulation has recently garnered great attention. Quantum Magnus algorithms, designed to achieve commutator scaling for time-dependent Hamiltonian simulation, have been found to be particularly efficient for such applications. When applied to unbounded Hamiltonian simulation in the interaction picture, they exhibit an unexpected superconvergence phenomenon. However, existing proofs are limited to the spatially continuous setting and do not extend to discrete spatial discretizations.
In this work, we provide the first superconvergence estimate in the fully discrete setting with a finite number of spatial discretization points $N$, and show that it holds with an error constant uniform in $N$. The proof is based on the two-parameter symbol class, which, to our knowledge, is applied for the first time in algorithm analysis. The key idea is to establish a semiclassical framework by identifying two parameters through the discretization number and the time step size rescaled by the operator norm, such that the semiclassical uniformity guarantees the uniformity of both. This approach may have broader applications in numerical analysis beyond the specific context of this work.\blfootnote{Emails: yonah\_borns-weil@berkeley.edu; di.fang@duke.edu; jiaqi.zhang988@duke.edu.} 
\blfootnote{ The manuscript has been accepted by \textit{Communications in Mathematical Physics}.}
\end{abstract}

\tableofcontents

\section{Introduction}

Hamiltonian simulation, a fundamental task in quantum computing, involves simulating the time evolution of a quantum system governed by a Hamiltonian. Originating as a motivation for the development of quantum computers, Hamiltonian simulation has remained central due to its crucial applications in physics and chemistry and its role as a subroutine in various quantum algorithms. The goal is to approximate the evolution described by the time-dependent Schrödinger equation:
\begin{equation}\label{eqn:ham_sim}
i\partial_t |\psi(t)\rangle = H(t) |\psi(t)\rangle, \quad |\psi(0)\rangle = |\psi_0\rangle,
\end{equation}
where $H(t)$ is the Hamiltonian. The exact evolution operator of this system is given by 
\begin{equation}  
U(t, 0) = \mathcal{T} \exp\left(-i \int_0^t H(s) \, ds\right),
\end{equation}
with $ \mathcal{T} $ denoting the time-ordering operator. When $ H(t) \equiv H $ is time-independent, this reduces to $U(t, 0) = \exp(-iHt)$.
Efficient Hamiltonian simulation for unbounded operators poses significant computational challenges, as the operator norm is infinite, resulting in substantial cost overhead under the typical operator norm error metric even after suitable discretizations.

In recent years, significant advances have been made in addressing these challenges of quantum algorithms. One notable approach involves Trotter-related algorithms, such as the Trotter product formulas, which decompose the unitary evolution operator into sequential exponentials of simpler terms, and multi-product formulas, which employ linear combinations of Trotter circuits. These approaches are known to achieve commutator scaling for time-independent Hamiltonians~\cite{ChildsSuTranEtAl2020,ChildsSu2019,WatkinsWiebeRoggeroLee2022,ZhukRobertsonBravyi2023,AftabAnTrivisa2024,Watson2025,WatsonWatkins2024},
delivering exceptional performance in applications to many quantum Hamiltonians.
Generalized Trotter formulas (or generalized Trotter-Suzuki expansion), which break the unitary evolution into simpler components of time-dependent Hamiltonian simulation \REV{(with the time-ordering still in their expressions and without further implementing such time-ordering operators yet) can} exhibit commutator scaling. This was first demonstrated for low-order cases~\cite{HuyghebaertDeRaedt1990,AnFangLin2021,RajputRoggeroWiebe2021} and later extended to arbitrarily high orders~\cite{MizutaTatsuhikoKeisuke2024}. \REV{Interestingly, the error bounds in high-order generalized Trotter formulas contain polynomial dependence on the time derivative of the Hamiltonian, not necessarily in a commutator form. Therefore, these high-order formulas benefit from commutator scaling when the time variation of the Hamiltonian is negligible, and are thus not ideal candidates for use in the interaction picture as we consider here. Additionally,}
each time-dependent building block still requires the inclusion of the time-ordering operator unless every sub-term in the Hamiltonian is in some specific form such as a linear combination of simple products of scalar time-dependent functions and time-independent matrices. Moreover, it has been shown that \REV{further} implementing each time-dependent building block \REV{in the generalized Trotter formulas} can result in the loss of the commutator scaling~\cite{AnFangLin2021}. 
On the other hand, post-Trotter methods~\cite{BerryChildsCleveEtAl2015,KieferovaSchererBerry2019,LowWiebe2019,BerryChildsSuEtAl2020,LowChuang2017,LowChuang2019,GilyenSuLowEtAl2019,ZlokapaSomma2024}, such as truncated Dyson series, quantum signal processing, and quantum singular value transformation, can achieve (near-)optimal scaling for bounded Hamiltonians. However, they still incur polynomial cost dependence on the operator norm due to the block-encoding subnormalization factor.

These challenges underscore the need for techniques that reduce dependence on operator norms, especially for the efficient simulation of unbounded Hamiltonians, which arise naturally in applications such as molecular and electronic calculations~~\cite{KivlichanWiebeBabbushEtAl2017,Somma2015,KivlichanMcCleanWiebeEtAl2018,AnFangLin2021,SuBerryWiebeEtAl2021,ChildsLengEtAl2022,RubinBerryKononovMaloneEtAl2023}, harmonic oscillators~\cite{Somma2015}, quantum optimization solvers~\cite{ZhangLengLi2021,LiuSuLi2023,LengHickmanLiWu2023,LengZhangWu2023,AugustinoLengNanniciniTerlakyWu2023}, and quantum field theories simulation~\cite{Preskill2019,TongAlbertMccleanPreskillSu2022,AbrahamsenTongBaoSuWiebe2023,PengSuClaudinoKowalskiLowRoetteler2023,SpagnoliWiebe2024,KuwaharaVuSaito2024}. Current analyses often focus on leveraging specific-case problem structures, e.g., in the observables, quantum states, or the operator algebras~\cite{Somma2015,SahinogluSomma2020,AnFangLin2021,SuHuangCampbell2021,ZhaoZhouShawEtAk2021,ChildsLengEtAl2022,FangTres2023,BornsWeilFang2022,HuangTongFangSu2023,ZengSunJiangZhao2022,GongZhouLi2023,LowSuTongTran2023,ZhaoZhouChilds2024,YuXuZhao2024,ChenXuZhaoYuan2024}).
Unbounded Hamiltonians have an infinite operator norm, posing computational challenges even when spatial discretizations are applied.
For instance, when discretizing the problem into a finite-dimensional space for numerical treatment, the matrix size and the operator norm of the discretized Hamiltonian both scale polynomially with respect to $N$, where $N$ represents the degrees of freedom in spatial discretization. 
A polynomial cost dependence on the operator norm would significantly deviate from the desired polylogarithmic complexity. 
This challenge is evident in finite-difference discretization examples, as shown in \cref{eq:A_fd}. 

Motivated by the above applications, we consider the quantum simulation of the Schr\"odinger equation, that is, the following unbounded Hamiltonian simulation problem:
\begin{equation}\label{eq:schd}
    i\partial_t \ket{\psi(t)} = H \ket{\psi(t)}, \quad H = -\frac{1}{2}\Delta + V(x),
\end{equation}
where $\Delta$ is the Laplacian operator, $V(x)$ is the potential, $x \in \mathbb{R}^d$ and $d$ is the dimension.
The simulation of the Schr\"odinger operator is BQP-hard~\cite{ZhengLengLiuWu2024}, and another interesting low\REV{er} bound considers the no-fast-forwarding theorem for a class of unbounded operators~\cite{TongAlbertMccleanPreskillSu2022}. When $V(x)$ has a bounded $L^\infty$ norm, the polynomial dependence on $N$ of quantum algorithms can often be resolved by considering the interaction picture~\cite{LowWiebe2019}, as discussed in \cite{AnFangLin2022,ChildsLengEtAl2022,FangLiuSarkar2024}. It turns unbounded Hamiltonian simulation into a time-dependent Hamiltonian simulation problem with the Hamiltonian
\begin{equation}
H_I(t) = e^{iAt} B e^{-iAt},
\end{equation} 
where $A$ and $B$ are numerical discretizations of the Laplacian $-\frac{1}{2}\Delta$ and the potential $V$, respectively.
One interesting line of work based on the quantum Magnus expansion~\cite{AnFangLin2022,FangLiuSarkar2024} -- which has been shown to exhibit commutator scaling for any time-dependent Hamiltonians -- has revealed a surprising phenomenon known as superconvergence. Other significant work along the line of quantum algorithms based on the Magnus series mainly for bounded Hamiltonians include \cite{SharmaTran2024,BosseChildsEtAl2024,CasaresZiniArrazola2024}. 

In numerical analysis, superconvergence refers to an algorithm performing significantly better than the expected accuracy in specific scenarios. In particular, a carefully designed quantum algorithm based on the second-order Magnus expansion can achieve fourth-order accuracy, with an error constant independent of the spatial discretization degrees of freedom $N$ for smooth potentials. This results in a computational cost that depends only polylogarithmically on $N$, while preserving the fourth-order accuracy. This efficiency is attributed to the time-dependent commutator scaling inherent in the algorithm for general time-dependent simulation, which facilitates additional error cancellation. Another impressive line of work, distinct from the concept of superconvergence yet closely related, is the study of destructive interference~\cite{TranChuSuEtAl2020, Layden2022}, which primarily focuses on bounded Hamiltonians, such as local lattice Hamiltonians.

In our previous study~\cite{FangLiuSarkar2024}, we showed that 
the second-order algorithm has superconvergence of fourth order. 
The analysis has treated the spatial variables as continuous, which is typically sufficient for spectral discretizations in numerical analysis. This also ensures in the continuum limit (as $N \to \infty$), the error remains bounded, and hence the cost remains insensitive to the increase of $N$. However, this does not guarantee that the behavior remains uniform with respect to $N$ for finite $N$ in discretizations such as finite differences. It turns out to treat finite $N$, previous techniques based on the pseudodifferential operators and semiclassical analysis in terms of simple symbol classes \REV{(see~\cref{eq:def_S(1)} for the simple symbol class definition)} no longer apply.

\bigskip
\noindent \textbf{Challenges and Contributions:}  

As discussed, previous results primarily apply to continuous space variables, making them suitable for spectral methods. A natural question arises: 
\begin{center}
    \textit{
    Can the superconvergence result be established for finite $N$ when working directly with discretized matrices in discretization methods such as the finite difference?
} 
\end{center}
This is the open question that we are about to address in this work. Addressing this question is in fact far from straightforward. \REV{At a high level, prior analyses of superconvergence for unbounded Hamiltonians crucially rely on a proof strategy based on the exact Egorov theorem, which no longer holds in the discrete case. In the following paragraphs, we explain in more detail what this means, why this issue is fundamental, and why it does not admit simple fixes.}

The core reason why extending these results to finite difference discretization is non-trivial lies in the exactness of Egorov’s theorem, which plays a crucial role in continuous-variable analyses. This theorem leverages properties intrinsic to the underlying (spatially continuous) operators that can be treated as pseudo-differential operators. For instance, when dealing with continuous operators like the Laplacian $A = -\frac{1}{2}\Delta$ and potential $B = V(x)$, Egorov’s theorem ensures that the following expression is exact (!), without truncation errors as seen in series expansions:
\begin{equation}
        e^{ -i\frac{s}{2}\Delta} V e^{i\frac{s}{2} \Delta} = e^{-i\frac{1}{2} th_0 \Delta} V e^{i \frac{1}{2}th_0 \Delta} = \mathrm{Op}_{h_0}(V(x + pt)), \quad s = th_0, \quad t\in[-1, 1].
\end{equation}
where $\mathrm{Op}_{h_0}(\cdot)$ is the Weyl quantization, which turns a phase-space function into the corresponding quantized operator (see the definition in \cref{{eq:weylquant}}), and $h_0 \in (0, 1]$ is an absolute constant that can be chosen as 1 in the proofs (see \cite{FangLiuSarkar2024}). This is essential in showing the error depends only on the potential $V$.

However, this exactness breaks down when transitioning to finite difference discretization. The operator symbol shifts from a quadratic form, such as $p^2$, to an approximation as $1 - \cos(p)$. While this form behaves like $p^2$ near $p = 0$, it deviates significantly from a true quadratic form for larger $p$, undermining the validity of the exact Egorov’s theorem. The exact version of the theorem is crucial for avoiding unboundedness in error estimates in the superconvergence case but does not hold in our discrete settings. 
Moreover, while an approximate version of Egorov’s theorem~\cite{zworski2022semiclassical,Martinez2002,Hormander2007} does exist, it holds only up to the Ehrenfest time, which scales as $\Or(\log(1/h)) = \Or( \log(N) )$. This constraint proves insufficient for the discrete case, which demands significantly longer time scales of $\Or(h^{-2}) = \Or(N^2) $.

To address the new challenges, two key questions must be considered:
First, are there discrete counterparts for pseudodifferential operators and semiclassical symbols that enable direct work with discrete matrices?
Second, if such counterparts exist, how can one address the loss of Egorov’s theorem and account for long-time scaling behavior to establish the superconvergence estimate? It turns out the first question can be answered via discrete microlocal analysis, introduced to numerical analysis and quantum algorithm estimates in \cite{BornsWeilFang2022}. This framework maps discrete matrices to the Weyl quantization of continuous symbols (phase-space functions) on the quantized torus. It has been shown that the matrix operator norms are bounded above by the operator norms of these continuous microlocal operators obtained by quantizing the phase-space functions~\cite{BornsWeilFang2022, DyatlovJezequel2021}. This approach effectively enables the use of (continuous) microlocal and PDE tools to directly analyze discrete matrices in numerical discretizations.

The second question, addressing the breakdown of Egorov's theorem, is far more critical and represents the most challenging aspect. To overcome it, we develop a new approach based on a two-parameter semiclassical symbol~\cite{NonnenmacherSjostrandZworski2014}. The key idea is to transform the problem through a unitary transformation and then show that the operators are in fact in the $\tilde{S}_{1/2}$ two-parameter symbol class \REV{(see its definition in~\cref{S_1/2})}. We then identify the inverse of matrix size $1/N$ as a semiclassical parameter (denoted as $h$), and incorporate the characteristic time scaling informed by the matrix operator norm as a secondary semiclassical parameter. This step was crucial as it enabled us to apply semiclassical calculus pertaining to the $\tilde{S}_{1/2}$ symbol class, known for providing estimates that are uniform in both semiclassical parameters. In our context, uniformity in $h$ translates directly to uniformity in $N$, which is precisely what we need for the fully discrete superconvergence estimate. \REV{We present in~\cref{sec:necessity_two_parameter_symbol} the necessity for the two-parameter symbol class and the high-level intuition of the proof argument.}

It is also worth noting that, unlike traditional settings where semiclassical analysis is naturally applicable due to the presence of an inherent parameter $h$ in the differential equation itself, our original problem does not include such a parameter.
By constructing a semiclassical framework and linking it to the discrete matrix problem, we were able to leverage uniform semiclassical estimates in $h$ to achieve uniform estimates in $N$.
This novel approach of connecting semiclassical results to discrete numerical algorithms may hold independent value for broader interests, providing a new framework for addressing similar challenges in the error estimates with spatial discretizations.

As a byproduct, we provide a much cleaner and conceptually simpler proof for the time-dependent commutator scaling in the quantum Magnus-based algorithm developed in \cite{FangLiuSarkar2024}. The previous proof relies on the remainder expansion of the Magnus series \cite[Lemma 2]{BlanesCasasOteoRos2009}, which is also commonly used for classical algorithms based on the Magnus series. This approach results in an estimation of 11 terms and, in addition to pseudodifferential arguments, some terms necessitate further Fourier analysis treatment. In contrast, our new analysis provides a more streamlined and efficient proof.
We derive an exact error representation in a compact form, where all terms naturally appear in the commutator structure $[H(t), [H(s), H(\tau)]]$ for any time-dependent Hamiltonian $H(t)$, without involving any remainder terms. This exact error representation of the Magnus series may also be of independent interest.

\bigskip
\noindent \textbf{Organization:}  

The rest of this paper is now organized as follows: In \cref{sect:algo_prelim}, we introduce the algorithmic preliminaries, including the interaction picture formulation and finite difference discretization. This section provides the necessary background on the spatial discretization of the differential operators and revisits the quantum Magnus algorithms for time-dependent Hamiltonian simulation. In \cref{sect:exact_error}, we establish an exact error representation for the second-order Magnus expansion.
We also present a derivation of local and long-time error bounds, in terms of the nested commutators.
\cref{sect:superconv} rigorously establishes superconvergence in the discrete setting, addressing the challenges posed by finite difference discretization and the breakdown of the Egorov’s theorem. We employ discrete microlocal analysis which turns the discrete superconvergence result into a continuous symbol estimate, which is then proved in \cref{sect:symbol_calculus}. 
\cref{sect:symbol_calculus} also introduces the two-scale symbol class, which are key mathematical tools to prove our result.
Finally, in \cref{sect:conclusion}, we conclude with some further remarks.

\section{Algorithmic Preliminaries}\label{sec:algorithm_prelim} \label{sect:algo_prelim}

In this section, we revisit the algorithmic preliminaries for the superconvergence of unbounded Hamiltonian simulation in the interaction picture. We first discuss the interaction picture, especially for the finite difference discretization of the Schr\"odinger equation, and then revisit the quantum Magnus algorithms. 

\subsection{Interaction Picture and Finite Difference Discretization}
The idea of the interaction picture can be applied to manage the unbounded-plus-bounded Hamiltonians, allowing the effective norm to remain finite. 

Specifically, when the Hamiltonian $H$ can be decomposed as $H = A + B$ with $\|A\| \gg \|B\|$, the interaction picture Hamiltonian~\cite{LowWiebe2019} is defined as:
\begin{equation}
H_I(t) = e^{iAt} Be^{-iAt},
\end{equation}
where $A$ represents the larger component and $B$ is the part with a smaller norm. This transformation ensures that $\|H_I(t)\| = \|B\| \ll \|A\|$, facilitating simulations with reduced complexity. In general, the interaction picture formulation also works for the Hamiltonian in the form of $f(t) A + B(t)$, where $f$ is a scalar function and $A$ is fast-forwardable.

When applying to the simulation of the Schr\"odinger equation~\cref{eq:schd}, $A$ and $B$ represent the discretization of the Laplacian and potential operators, respectively. In our work, we consider the finite difference discretization.  Note that spectral discretization provides spectral accuracy, and the spatially continuous proofs presented in \cite{AnFangLin2022, FangLiuSarkar2024} are already sufficient. \REV{We note that, typically in numerical analysis, when considering spectral spatial discretizations, it is acceptable to work in the semi-discrete setting (i.e., only discretizing in time). 
This is because spectral methods achieve spectral accuracy and, from the perspective of numerical analysis, preserve properties such as the product rule and chain rule, similarly to continuous operators--unlike finite difference methods, where the product rule does not hold in the same way as in the continuous case.
} However, for finite difference discretization, the discrete matrix might exhibit behavior much more distinct from its continuous counterpart. As will be made clear later in \cref{sec:finite_difference_discretization}, the corresponding symbol undergoes significant changes compared to its continuous counterpart, leading to the failure of previous arguments and an entirely new proving mechanism.

We now focus on the discretization of the Schr\"odinger equation by using a finite difference method with $N$ evenly spaced nodes $x_j = a + (b - a)j/N $ where $0 \leq j \leq N-1$ with periodic boundary conditions. Note that periodic boundary conditions are commonly employed in the simulation of the Schrödinger equation, as periodization can often be applied to reformulate the problem in this format (see, for example, the textbook~\cite{Lubich2008book}).
This gives the matrices $ A$ and $B$ that represent the discretized versions of the operators $ -\frac{1}{2}\Delta$ and $V(x)$ in one dimension, given by:
\begin{equation}\label{eq:A_fd}
   A =  \tilde{A}_{\mathrm{dis}} :=  \frac{N^2}{2(b-a)^2} \left(\begin{array}{ccccc}
        2 & -1 & & & -1 \\
         -1& 2 & -1 & & \\
          & \ddots& \ddots& \ddots& \\
           & & -1& 2 & -1\\
        -1& & & -1 & 2 \\
    \end{array}\right)_{N \times N},
\end{equation}
and 
\begin{equation}\label{eq:B_fd}
    B = \text{diag}\left(V(x_0),V(x_1),\cdots,V(x_{N\REV{-1}}) \right).
\end{equation}
In $d$ dimension, one uses $N$ grids for each dimension. $B$ remains a diagonal matrix, while $A$ becomes 
\begin{equation}\label{eq:A_fd_hd}
A= \tilde{A}_{\mathrm{dis}}\otimes I_N^{\otimes d-1}+ I_N \otimes \tilde{A}_{\mathrm{dis}}\otimes I_N^{\otimes d-2} \dots +I_N^{\otimes d-1}\otimes \tilde{A}_{\mathrm{dis}},
\end{equation}
where $I_N$ is the identity matrix over $N$ dimension. \REV{Here, we focus on periodic boundary conditions in our work. In Schr\"odinger PDE simulations, it is common to apply the so-called periodization and then impose periodic boundary conditions (we refer the reader to the book \cite{Lubich2008book} for further discussion).}

Recent work has demonstrated that quantum algorithms based on truncating the Magnus series in this context can yield a surprising property known as superconvergence~\cite{AnFangLin2022,FangLiuSarkar2024}. This property is characterized by achieving higher convergence orders than expected (e.g., second-order truncation yielding fourth-order convergence). It holds with error preconstants that depend only on the potential $V$. This feature is particularly beneficial for real-space simulations, where maintaining computational efficiency is crucial. In this work, we will show that with finite $N$, such superconvergence estimates remain valid and uniform in $N$. 

\REV{Before discussing the quantum Magnus algorithm, we stress our motivation in considering the quantum Magnus algorithms. While there are other types of algorithms that can handle the interaction picture simulation, notably the truncated Dyson series (see, e.g., \cite{BerryChildsCleveEtAl2015,LowWiebe2019,KieferovaSchererBerry2019,BerryChildsSuEtAl2020}). Although the truncated Dyson series does not introduce the polynomial dependence on the time derivative of the Hamiltonian, it does not have commutator scaling, and it does not exhibit the superconvergence phenomenon. So are other time-dependent simulation algorithms, such as the randomized algorithms and the classical Magnus algorithms (see, e.g.,~\cite{BerryChildsSuEtAl2020,PoulinQarrySommaVerstraete2011,BlanesCasasRos2000,BlanesCasasOteoRos2009,BlanesCasas2025}). At a high level, the superconvergence arises from additional cancellation of the error in a specific way for unbounded Hamiltonian simulation. Ultimately, such cancellation is made possible thanks to the commutator scaling property of the quantum Magnus algorithms. Indeed, among various types of time-dependent Hamiltonian simulation algorithms, the quantum Magnus algorithm is, to the best of our knowledge, the only one that exhibits commutator scaling while avoiding the reintroduction of polynomial dependence on the derivative for general time-dependent Hamiltonians. Therefore, it is currently the only algorithm expected to display such a surprising superconvergence behavior for unbounded operators in the interaction picture.
}

\subsection{Revisit of the quantum Magnus algorithm}\label{sec:revisit_quantum_algorithm}

In this section, we review the quantum algorithm based on the Magnus expansion for general time-dependent Hamiltonian simulation as proposed in \cite{FangLiuSarkar2024}. 
This quantum algorithm is derived from the Magnus series, which approximates a time-order exponential via the exponential of a series without the time ordering.
Specifically, the solution of the linear differential equation with time-dependent coefficients $\REV{G}(t)$ (for Hamiltonian simulation $\REV{G}(t) = -i H(t)$)
\begin{equation}
    i\partial_t u(t) = H(t) u(t), \quad u(0) = u_0,
\end{equation}
can be represented as
\begin{equation}
  u(t)=\exp(\Omega (t))u_0
  , \quad  
   \Omega(t) =\sum_{n=1}^\infty \Omega^{(n)}(t),
\end{equation}
where the term in the Magnus series then reads
\begin{equation}
    \Omega^{(1)}(t) = \int_0^t \REV{G}(t_1)dt_1
\end{equation}
\begin{equation}   \label{eqn:mag_n}
  \Omega^{(n)}(t) =  \sum_{j=1}^{n-1} \frac{B_j}{j!} \,
    \sum_{
            k_1 + \cdots + k_j = n-1 \atop
            k_1 \ge 1, \ldots, k_j \ge 1}
            \, \int_0^t \,
       \ad_{\Omega^{(k_1)}(s)} \,  \ad_{\Omega^{(k_2)}(s)} \cdots
          \, \ad_{\Omega^{(k_j)}(s)} \REV{G}(s) \, ds    \qquad n \ge 2.
\end{equation}
Here $B_n$ are the Bernoulli numbers and $B_1 = -1/2$, and the adjoint action $\ad$ is defined via the commutators $\ad_{\Omega}(C) = [\Omega, C]$ and $\ad^{k+1}_{\Omega}(C) = [\Omega, \ad^{k}_{\Omega}(C)]$ for $k \in \mathbb{N}$. For an in-depth exploration of the Magnus series, readers can refer to the comprehensive review in, e.g., ~\cite{BlanesCasasOteoRos2009}. For notational simplicity, we denote the sum of the first $n$ terms as $\Omega_n$, namely,
\begin{equation}
    \Omega_n = \sum_{j=1}^n \Omega^{(j)}, \quad \Omega (t) = \lim_{n\to \infty}  \Omega_n(t).
\end{equation}
It is convenient to observe that the few terms of $\Omega$ are
\vspace{-.5em}
 \begin{align}\label{eq:omega_1_omega_2}
  \Omega^{(1)}(t) &= \int_0^t \REV{G}(t_1)\,dt_1, \quad
  \Omega^{(2)}(t) = \frac{1}{2} \int_0^t dt_1 \int_0^{t_1} dt_2 \, [\REV{G}(t_1), \REV{G}(t_2)], 
\end{align}

The quantum algorithm can be then summarized as follows. We divide the time interval of $[0, T]$ into $L$ equidistant time sub-intervals.
Within each short time interval $[t_j, t_{j+1}]$ (where $\Delta t = T/L$ represents the time step and $t_j = j \Delta t$), the time-ordered matrix exponential is approximated using the second-order Magnus expansion $\Omega_2(t_{j+1}, t_j)$,
\begin{equation}\label{eq:def_Uexact_U2_tj}
    U_\mathrm{exact}(t_{j+1},t_j) = \mathcal{T} e^{-i \int_{t_j}^{t_{j+1}} H(s) ds} \approx e^{\Omega_2(t_{j+1}, t_j)} : = U_2(t_{j+1}, t_j),
\end{equation}
where
\begin{equation}\label{eq:def_omega2_tj}
    \Omega_2(t_{j+1}, t_j): = -i\int_{t_j}^{t_{j+1}} H(s) ds + \frac{1}{2} \int_{t_j}^{t_{j+1}} \left[
  \int_{t_j}^{s} H(\sigma) d\sigma, H(s) \right] ds.
\end{equation}
 The integrals are implemented accurately and efficiently via a quantum circuit based on the linear combination of unitaries (LCU)~\cite{ChildsWiebe2012,GilyenSuLowEtAl2019} and compare oracle~\cite{Oliveira2007QuantumBit,sanders2019black}. The idea is to conduct the numerical quadrature via Riemann sums with a large number of quadrature points $M$. In particular, we need $M$ to be polynomially dependent on $\|\partial_ t H(t)\|$. However, the quantum implementation via the circuits \cite[Section 5.2 and Fig. 2]{FangLiuSarkar2024} only requires a cost (the number of quantum gates) that is $\log(M)$ and hence a cost only logarithmically dependent on the time derivative of $H(t)$. This distinguishes itself from many other Magnus series based numerical integrators, especially classical ones~\cite{HochbruckLubich2003,Thalhammer2006,BlanesCasasThalhammer2017}, which have a cost polynomially dependent on the time derivative $H(t)$.

 This feature is crucially important in the context of unbounded Hamiltonian simulation considered here. Throughout the work, we consider the time-dependent Hamiltonian as the Hamiltonian in the interaction picture, given by
 \begin{equation}\label{eq:Ham_interaction_picture}
     H(t) = H_I(t) =  e^{iAt} Be^{-iAt},
 \end{equation}
 where $A$ and $B$ are finite difference discretizations of the Laplacian and potential operators. From this, we observe that
 \begin{equation}
     \partial_t H_I(t)  = i e^{iAt} [A, B]e^{-iAt},
 \end{equation}
implying
 \begin{equation}
     \norm{\partial_t H_I(t) } = \norm{[A, B]} = \Or(N).
 \end{equation}
 If an algorithm’s cost depends polynomially on 
$\norm{\partial_t H_I(t) }$, the overall cost would scale polynomially with $N$. However, achieving a polylogarithmic dependence on $N$ is the primary focus here. Therefore, it is essential that time-dependent Hamiltonian simulation algorithms maintain a polylogarithmic cost in the time derivative of $H(t)$ and hence in $N$, which is a nice feature for the quantum algorithm we study here. Thanks to this feature of the quantum implementation, proving $N$-independent error bounds for the difference between $U_\mathrm{exact}$ and $U_2$ is sufficient to establish the superconvergence result. Therefore, we focus on estimating the difference between $U_\mathrm{exact}$ and $U_2$ in the following analysis.

\section{Exact Error Representation} \label{sect:exact_error}

Using the variation of constants formula, we can explicitly write out the exact error representation for the case of the second-order Magnus expansion.
On one hand, the exact unitary $U(t)$ satisfies
\begin{equation}
 \partial_t U_\mathrm{exact}(t) =\REV{G}(t)U_\mathrm{exact}(t).
\end{equation}
On the other hand, recall that the second-order Magnus expansion $U_2(t)$ is given as
\begin{equation*}
 U_2(t)=e^{\int_0^t \REV{G}(s)\,ds-\frac{1}{2}\int_0^t ds \int_0^s d\sigma[\REV{G}(\sigma),\REV{G}(s)]}=e^{\Omega_2(t)}. 
\end{equation*}
To represent the difference between $U_\mathrm{exact}(t)$ and $U_2(t)$, we consider the differential equation satisfied by $U_2$, which we denote as
\begin{equation}
  \partial_t {U}_2=\REV{\tilde{G}}(t)U_2(t).
\end{equation}
An expression of $\REV{\tilde{G}}$ will subsequently provide the error representation through another application of the variation of constants formula.

\REV{We remark that although the second-order Magnus error analysis can be done in various different ways, such as utilizing a Taylor theorem approach (with remainders) to represent this operator $\REV{\tilde{G}}(t)$, as presented in \cite{FangLiuSarkar2024}. The resulting representation would involve a remainder term in the Taylor series, namely,
\begin{equation}\label{eqn:tilde_A}
    \REV{\tilde{G}}(t) =  d \exp_{\Omega_2(t)}(\dot \Omega_2) = \dot \Omega_2 + \frac{1}{2}[\Omega_2, \dot \Omega_2] + \frac{1}{6}[\Omega_2, [\Omega_2, \dot \Omega_2]] +g(\ad_{\Omega_2})(\ad_{\Omega_2}^3(\dot \Omega_2) ),
\end{equation}
where the last term is given by the remainder in the series expansion of the function 
\begin{equation}
    \frac{e^{z}- 1}{z} = 1 + \frac{1}{2}z + \frac{1}{6} z^2 +  g(z) z^{3}.
\end{equation}
While for bounded operators everything is rigorous and well-defined, for unbounded operators one also needs to check that the operator $g(\ad_{\Omega_2})(\ad_{\Omega_2}^3(\dot{\Omega}_2))$ is defined from the correct domain to the correct range, which was carefully proved in \cite[Appendix C]{FangLiuSarkar2024}. For the discrete superconvergence case, although all operators are discretized and hence bounded, we are interested in their uniformity with respect to $N$. One also needs to perform the uniform-in-$N$ estimate for this remainder operator governed by $g$, which--while possible--can be rather complicated (and, in fact, unnecessary). Here, we present a different error representation strategy. As can be seen in \cref{eq: G(t)}, we write an integral representation of $\REV{\tilde{G}}$ instead of using the previous Taylor expansion approach, which does not involve the term $g(\ad_{\Omega_2})(\ad_{\Omega_2}^3(\dot{\Omega}_2))$.
  }

First, recall the following elementary lemma for the derivative formula of the matrix exponential without time ordering: 
\begin{lemma}[Derivative of Matrix Exponential without time ordering] \label{lem:derivative_matrix_exponential}
\begin{equation}
\begin{aligned}
\frac{\mathrm{d}}{\mathrm{d}t} e^{X(t)} = \int_0^1 e^{\lambda X(t)} \frac{\mathrm{d}X}{\mathrm{d}t} e^{(1-\lambda)X(t)} \, \mathrm{d}\lambda.
\end{aligned}
\end{equation}
\end{lemma}
This formula is well-known, and we do not include the proof here. It can be easily derived using the variation of constants formula. \REV{To apply the lemma in our case, we need an expression of the time derivative of $\Omega_2$, namely,}
\begin{equation*}
   \frac{\mathrm{d}}{\mathrm{d}t} \Omega_2=\REV{G}(t)-\frac{1}{2}\int_0^t d\sigma[\REV{G}(\sigma), \REV{G}(t)],
   \end{equation*}
\REV{which follows from \cref{eq:omega_1_omega_2}.}
 Applying \cref{lem:derivative_matrix_exponential}, it follows that
 \begin{equation*}
 \frac{\mathrm{d}}{\mathrm{d}t} U_2 = \frac{\mathrm{d}}{\mathrm{d}t} e^{\Omega_2(t)}= \int_0^1 \mathrm{d}\lambda \; e^{\lambda \Omega_2(t)} \Bigg(\REV{G}(t) - \frac{1}{2} \int_0^t \mathrm{d}\sigma [\REV{G}(\sigma), \REV{G}(t)] \Bigg)e^{-\lambda \Omega_2(t)} U_2.
 \end{equation*}
 \begin{equation*}
    \REV{\tilde{G}}(t) - \REV{G}(t) = \int_0^1 \mathrm{d}\lambda \; e^{\lambda \Omega_2(t)}\REV{G}(t)e^{-\lambda \Omega_2(t)} - \REV{G}(t)-\frac{1}{2} e^{\lambda \Omega_2(t)}\int_0^t \mathrm{d}\sigma [\REV{G}(\sigma), \REV{G}(t)] e^{-\lambda \Omega_2(t)}. 
\end{equation*}
In particular,
\begin{equation*}
    e^{\lambda \Omega_2(t)}\REV{G}(t)e^{-\lambda \Omega_2(t)}-\REV{G}(t)=\int_0^{\lambda}\, d\tilde{\lambda}\frac{d}{d\tilde{\lambda}}\left( e^{\tilde{\lambda} \Omega_2(t)} \REV{G}(t) e^{-\tilde{\lambda}\Omega_2(t)}\right)=\int_0^{\lambda}e^{\tilde{\lambda} \Omega_2(t)}[\Omega_2(t),\REV{G}(t)]e^{-\tilde{\lambda} \Omega_2(t)}\,d{\tilde{\lambda}}.
\end{equation*}
Therefore,
\begin{equation}\label{eq: G(t)}
\begin{aligned}
\REV{\tilde{G}}(t) - \REV{G}(t) &= \int_0^1 d\lambda \int_0^\lambda d\tilde{\lambda} \, e^{\tilde{\lambda} \Omega_2(t)} [\Omega_2(t), \REV{G}(t)] e^{-\tilde{\lambda} \Omega_2(t)} \\
&- \frac{1}{2} \int_0^1 d\lambda \, e^{\lambda \Omega_2(t)} \int_0^t d\sigma \, [\REV{G}(\sigma), \REV{G}(t)] e^{-\lambda \Omega_2(t)}.
\end{aligned}
\end{equation}
Since we know that
\begin{equation}\label{eq: commutator}[\Omega_2(t),\REV{G}(t)]=\int_0^t \,d\sigma [\REV{G}(\sigma),\REV{G}(t)]-\frac{1}{2}\int_0^t\,ds\int_0^s\,d\sigma\left[[\REV{G}(\sigma),\REV{G}(s)],\REV{G}(t)\right],
\end{equation}
After substituting \cref{eq: commutator} into \cref{eq: G(t)}, it follows that 
\begin{equation*}
\begin{aligned}
   &\int_{0}^{1} d\lambda \int_{0}^{\lambda} d\tilde{\lambda} \, e^{\tilde{\lambda} \Omega_{2}(t)}  \int_{0}^{t} d\sigma [\REV{G}(\sigma), \REV{G}(t)] e^{-\tilde{\lambda} \Omega_{2}(t)}\\
&- \frac{1}{2} \int_{0}^{1} d\lambda \, e^{\lambda \Omega_{2}(t)} 
\int_{0}^{t} d\sigma [\REV{G}(\sigma), \REV{G}(t)] e^{-\lambda \Omega_{2}(t)}\\
&-\frac{1}{2}\int_{0}^{1} d\lambda \int_{0}^{\lambda} d\tilde{\lambda} \, e^{\tilde{\lambda} \Omega_{2}(t)} \int_0^t\,ds\int_0^s\,d\sigma\left[[\REV{G}(\sigma),\REV{G}(s)],\REV{G}(t)\right] e^{-\lambda \Omega_{2}(t)}\\
&:=I_1+I_2+I_3.
 \end{aligned}
 \end{equation*}
  where two of the three terms in $\REV{\tilde{G}}(t) - \REV{G}(t)$ are 
 \begin{equation*}
     I_1=\int_{0}^{1} d\lambda \int_{0}^{\lambda} d\tilde{\lambda} \, e^{\tilde{\lambda} \Omega_{2}(t)}  \int_{0}^{t} d\sigma [\REV{G}(\sigma), \REV{G}(t)] e^{-\tilde{\lambda} \Omega_{2}(t)},
 \end{equation*}
 \begin{equation*}
 \begin{aligned}
     I_2= & -\frac{1}{2} \int_{0}^{1} d\lambda \, e^{\lambda \Omega_{2}(t)} 
\int_{0}^{t} d\sigma [\REV{G}(\sigma), \REV{G}(t)] e^{-\lambda \Omega_{2}(t)}
\\
= &-\frac{1}{2}\int_0^1\,d\lambda \int_0^1d\tilde{\lambda}\,e^{\lambda \Omega_{2}(t)} 
\int_{0}^{t} d\sigma [\REV{G}(\sigma), \REV{G}(t)] e^{-\lambda \Omega_{2}(t)}.
 \end{aligned}
 \end{equation*}
 Also note that $I_3$ is already a two-layer nested commutator with two layers of time integrals, as desired.
We define
\begin{equation}\label{eq:Bt}
    \int_{0}^{t} d\sigma [\REV{G}(\sigma), \REV{G}(t)]\triangleq \REV{K}(t),
\end{equation}
and
\begin{equation*}
    e^{\lambda \Omega_{2}(t)} 
\int_{0}^{t} d\sigma [\REV{G}(\sigma), \REV{G}(t)] e^{-\lambda \Omega_{2}(t)}\triangleq e^{\lambda \Omega_{2}(t)} \REV{K}(t) e^{-\lambda \Omega_{2}(t)}  \triangleq f(\lambda).
\end{equation*}
After simplification, since in the integral with respect to $\lambda$ and $\tilde{\lambda}$ is $f(\tilde{\lambda})-\frac{1}{2}f(\lambda)$, we have
\begin{equation}\label{eq: I1I2}
\begin{aligned}
I_1+I_2&=\int_0^1\,d\lambda\int_0^{\lambda} d\tilde{\lambda}f(\tilde{\lambda})-\frac{1}{2}\int_0^1\,d\lambda \,f(\lambda) \\
&=\int_0^1\,d\lambda\int_0^{\lambda} d\tilde{\lambda}\,f(\tilde{\lambda})-\frac{1}{2}\int_0^1\,d\lambda\int_0^1 d\tilde{\lambda}f(\lambda),
\end{aligned}
\end{equation}
from which we can see that $I_1+I_2$ is an integral of $f(\tilde{\lambda}) - \frac{1}{2}f(\lambda)$. We observe that if $f(\lambda)$ is independent of $\lambda$, say $f(\lambda) \equiv f$, then the integral would equal to 0, implying that there will be no error contribution from these terms. Therefore, only the dependence of $f(\lambda)$ on $\lambda$ contributes to $I_1+I_2$, and the derivation focuses entirely on the contributions from the higher-order nested commutators. To see this, \begin{equation*}
I_1+I_2=\int_0^1\,d\lambda\int_0^{\lambda} d\tilde{\lambda}f(\tilde{\lambda})-\frac{1}{2}\int_0^1\,d\lambda\int_0^1 d\tilde{\lambda}f(\lambda),    
\end{equation*} 
which implies that
\begin{equation*}
    I_1+I_2= \frac{\lambda^2}{2}f\,\bigg|_0^1 -\frac{1}{2} f=0.
\end{equation*}
In light of this, since we want the error terms to be all in the form of nested commutators, our goal now becomes analyzing whether $2f(\lambda)-f(\tilde\lambda)$ only contains terms with the nested commutator form of at least 3 layers. 
As defined before, the function $f(\lambda)$ is derived from the expression involving the operator $\REV{K}(t)$: $f(\lambda)=e^{\lambda\Omega_2(t)}\REV{K}(t)e^{-\lambda\Omega_2(t)}$, where $\Omega_2(t)$ is a commutator operator derived from earlier terms in the Magnus expansion.

We then perform the Taylor expansion of $f(\lambda)=e^{\lambda\Omega_2(t)}\REV{K}(t)e^{-\lambda\Omega_2(t)}$. Recall that 
\begin{equation}
\frac{d}{d\lambda} \left( e^{\lambda \Omega_2(t)} \REV{K}(t) e^{-\lambda \Omega_2(t)} \right) = e^{\lambda \Omega_2(t)} [\Omega_2(t), \REV{K}(t)] e^{-\lambda \Omega_2(t)}.
\end{equation}
Integrating both sides with respect to $\lambda$ gives us
\begin{equation}
\begin{aligned}
f(\lambda) &= \REV{K}(t) + \int_0^\lambda d\tilde{\lambda} \, e^{\tilde{\lambda} \Omega_2(t)} [\Omega_2(t), \REV{K}(t)] e^{-\tilde{\lambda} \Omega_2(t)}.
\end{aligned}
\end{equation}
Here $\REV{K}(t)$, defined in \cref{eq:Bt}, is independent of $\lambda$, so its contribution in $I_1+I_2$ is 0.
Thus, higher order terms (more layers nested commutator) show the dependency on the nested commutators. We thus keep only terms involving three or more layers of nested commutators for $\REV{\tilde{G}} - \REV{G}$, which are the integral remainder term for $I_1+I_2$ after Taylor expansion.
We proceed by substituting the Taylor series expansion of $f(\lambda)$ and $f(\tilde{\lambda})$ into the expression of $I_1+I_2$ as in \cref{eq: I1I2}:
\begin{equation}
\begin{aligned}
I_1+I_2 =& \int_0^1\,d\lambda\int_0^{\lambda} d\tilde{\lambda}\left(\REV{K}(t) + \int_0^{\tilde{\lambda}} d\mu \, e^{\mu \Omega_2(t)} [\Omega_2(t), \REV{K}(t)] e^{-\mu \Omega_2(t)} \right) \\
 &-\frac{1}{2}\int_0^1\,d\lambda\int_0^1 d\tilde{\lambda}\left(\REV{K}(t) + \int_0^{\lambda } d\mu \, e^{\mu \Omega_2(t)} [\Omega_2(t), \REV{K}(t)] e^{-\mu \Omega_2(t)} \right)\\
 =&\int_0^1 d\lambda \Bigg(  
 \int_0^{\lambda} d\tilde{\lambda} \int_0^{\tilde{\lambda}} d\mu \, e^{\mu \Omega_2(t)} [\Omega_2(t), \REV{K}(t)] e^{-\mu \Omega_2(t)}
 \\
 &  \qquad \qquad -\frac{1}{2}\int_0^1 \,d\tilde{\lambda}\int_0^{\lambda}  d\mu \, e^{\mu \Omega_2(t)} [\Omega_2(t), \REV{K}(t)] e^{-\mu \Omega_2(t)} 
 \Bigg).
 \end{aligned}
\end{equation}
As $\REV{K}(t)$ contains a one-layer commutator as defined in \cref{eq:Bt}, the term $[\Omega_2(t), \REV{K}(t)]$ contains a two-layer nested commutator. Hence, $I_1+I_2$ after the cancellation of $\REV{K}(t)$ gives us two layers of the time integrals with a two-layer nested commutator. In the later superconvergence estimate, we will show that the two-layer nested commutator in terms of $\REV{G}(t)$ contributes to an additional $\Or(t^2)$. 
As discussed before, $ I_3 $ is already in the form of two-layer nested commutators in terms of $\REV{G}(t)$. 
Combining all terms, one has:
    \begin{equation} \label{eq:tilde_A-A_final}
        \begin{aligned}    
   & \REV{\tilde{G}} - \REV{G} =  I_3 +I_1+I_2 
    \\
    =&  -\frac{1}{2}\int_{0}^{1} d\lambda \int_{0}^{\lambda} d\tilde{\lambda} \, e^{\tilde{\lambda} \Omega_{2}(t)} \int_0^t\,ds\int_0^s\,d\sigma\left[[\REV{G}(\sigma),\REV{G}(s)],\REV{G}(t)\right]e^{-\tilde{\lambda} \Omega_{2}(t)}
    \\
    &+ \int_0^1 d\lambda \Bigg(  
 \int_0^{\lambda} d\tilde{\lambda} \int_0^{\tilde{\lambda}} d\mu \, e^{\mu \Omega_2(t)} [\Omega_2(t), \REV{K}(t)] e^{-\mu \Omega_2(t)}
 \\
 & \qquad \qquad
 -\frac{1}{2}\int_0^1 \,d\tilde{\lambda}\int_0^{\lambda}  d\mu \, e^{\mu \Omega_2(t)} [\Omega_2(t), \REV{K}(t)] e^{-\mu \Omega_2(t)} 
 \Bigg).
    \end{aligned}
     \end{equation}
Recall that $[\Omega_2, \REV{K}(t)]$ is in a similar form as in \cref{eq: commutator}, namely,
\begin{equation}
\begin{aligned} 
& [\Omega_2(t),\REV{K}(t)]
\\ = & \int_0^t \,d\sigma [\REV{G}(\sigma),\REV{K}(t)]-\frac{1}{2}\int_0^t\,ds\int_0^s\,d\sigma\left[[\REV{G}(\sigma),\REV{G}(s)],\REV{K}(t)\right]
\\
=&  \int_0^t \,d\sigma \int_{0}^{t} d\tau [\REV{G}(\sigma),[\REV{G}(\tau), \REV{G}(t)]]-\frac{1}{2}\int_0^t\,ds\int_0^s\,d\sigma \int_{0}^{t} d\tau \left[[\REV{G}(\sigma),\REV{G}(s)],[\REV{G}(\tau), \REV{G}(t)]\right].
\end{aligned}
\end{equation}
The above calculation yields the exact error representation over the interval $[0,t]$. For the local truncation error on $[t_j, t_{j+1}]$, a similar analysis holds. One only needs to, for example, change the integral bounds from $[0,t]$ to $[t_j, t_{j+1}]$ or replace $\REV{G}(t)$ by $\REV{G}(t_j+t)$ with $t\in[0,\Delta t]$ to properly represent the value of $\REV{G}$ over $[t_j, t_{j+1}]$. Together with the following equation,
\begin{equation}
        U_\mathrm{exact}(t_{j+1}, t_j )  - U_2 (t_{j+1}, t_j)  = \int_{t_j}^{t_{j+1}} U_\mathrm{exact}(  t_{j+1}, s) 
    \left( \REV{G}(s) - \REV{\tilde{G}}(s) \right) U_2(s, t_j) \; ds,
\end{equation}
yields the exact error representation for the second-order Magnus expansion, which we do not expand further here. The key observation is that all terms are expressed in a nested commutator form with at least two layers, and this error representation is exact, obtained without using any Taylor series expansion that would introduce implicit remainder terms.

Now we proceed to find the upper bound of the norm of $U_\mathrm{exact}  - U_2$ as
   \begin{equation}
       \begin{aligned}
        &    \|U_\mathrm{exact}(t_{j+1}, t_j )  - U_2(t_{j+1}, t_j )\| \leq \Delta t \sup_{s \in [t_j,t_{j+1}]}\norm{ \REV{G}(s) - \REV{\tilde{G}}(s)}
           \\
           \leq \ &\frac{13}{24}\Delta t^3\sup_{\sigma, s, t\in [t_j,t_{j+1}]} \|\left[[\REV{G}(\sigma),\REV{G}(s)],\REV{G}(t)\right]\|
           \\
           & +\frac{5}{48}\Delta t^4 \sup_{\sigma, s, t, \tau\in [t_j,t_{j+1}]}\|  \left[[\REV{G}(\sigma),\REV{G}(s)],[\REV{G}(\tau), \REV{G}(t)]\right] \|,          %
       \end{aligned}
   \end{equation}
   where we also used the fact that the norms of unitary operators are 1.

We have reached the following theorem.

\begin{thm}[Local Truncation Error via Exact Error Representation] \label{thm:lte-mag2} 
Consider the exact unitary $U_\mathrm{exact}$, and the second-order Magnus expansion $U_2$ as defined in \cref{eq:def_Uexact_U2_tj}, \cref{eq:def_omega2_tj}. If for $t\in[t_j, t_{j+1}]$, the nested commutators satisfy 
\begin{equation}  \label{eq:nest_comm_alpha}
    \sup_{\sigma, s, t\in [t_j,t_{j+1}]} \|\left[[H(\sigma),H(s)],H(t)\right]\| \leq \alpha_\mathrm{comm}^j,
\end{equation}
\begin{equation}\label{eq:nest_comm_beta}
\sup_{\sigma, s, t, \tau\in [t_j,t_{j+1}]}\|  \left[[H(\sigma),H(s)],[H(\tau), H(t)]\right] \|\leq \beta_\mathrm{comm}^j,
\end{equation}
then we have 
\begin{equation}
    \norm{U_\mathrm{exact}( t_{j+1}, t_j) - U_2(t_{j+1}, t_j)} \leq \frac{13}{24}\Delta t^3 \alpha_\mathrm{comm}^j + \frac{5}{48}\Delta t^4 \beta_\mathrm{comm}^j. 
\end{equation}

\end{thm}

\begin{thm}[Long-time Error for Hamiltonian Simulation]\label{thm:long-time-mag2}  
Let 
\[
U_\mathrm{exact}(T,0) = \mathcal{T}e^{-i\int_0^T H(s) \, ds}
\]
be the exact long-time unitary over the time interval $[0, T]$, and the quantum Magnus algorithm accounting for the numerical quadrature $\Tilde{U}_2(T,0)$ be
\begin{equation}
    \Tilde{U}_2(t_L,t_{L-1})\dots \Tilde{U}_2(t_2,t_1)\Tilde{U}_2(t_1,t_0),
\end{equation}
where $\Tilde{U}_2(t_{j+1},t_j)=e^{\tilde{\Omega}_2(t_{j+1}, t_j)}$ and $\tilde{\Omega}_2$ represents the Riemann sum quadrature of $\Omega_2$ defined in \cref{eq:def_omega2_tj}. For $\norm{H(t)} \leq \alpha$ for all $t\in [0, T]$,
we have
\begin{equation}\label{eq:long_time_with_quadrature}
\begin{aligned}
    & \norm{U_\mathrm{exact}( T, 0) - \tilde{U}_2(T, 0)} 
    \\
    \leq & \frac{13}{24}\Delta t^3 \sum_{j=0}^{L-1}\alpha_\mathrm{comm}^j 
    + \frac{5}{48}  \Delta t^4\sum_{j=0}^{L-1}\beta_\mathrm{comm}^j
    +\frac{T \Delta t}{M} (1+3 \Delta t \alpha)\sup_{s\in [0,T]}\norm{H'(s)} ,
    \end{aligned}
\end{equation}
where $\alpha_\mathrm{comm}^j$ and $\beta_\mathrm{comm}^j$ are as defined in \cref{thm:lte-mag2}, and $M$ is the number of numerical quadrature points used in the Riemann sum.
\end{thm}

\begin{proof}

The proof follows standard argument in numerical analysis, due to each of $\Tilde{U}_2$ is unitary. Specifically,  let
\[
\| U(t_{j+1}, t_j) - \tilde{U}_j \| = \| r_j \| \leq e_j.
\]
We claim that for unitary 
$\tilde{U}_j$,
\begin{equation}
\| U(t, 0) - \tilde{U}_{L-1} \cdots \tilde{U}_0 \| \leq \sum_{j=0}^{L-1} e_j.
\end{equation}
This follows a straightforward induction
\begin{equation}
\begin{aligned}
&U(t_L, t_{L-1}) \cdots U(t_2, t_1)U(t_1, t_0)-\tilde{U}_{L-1}\cdots \tilde{U}_1\tilde{U}_0 
\\
=\ & U(t_L, t_{L-1}) \cdots U(t_2, t_1) \left( U(t_1, t_0)-\tilde{U}_0\right)
+\left(U(t_L, t_{L-1})\cdots U(t_2, t_1)-\tilde{U}_{L-1}\cdots \tilde{U}_1\right)\tilde{U}_0,
\end{aligned}
\end{equation}
which implies 
\begin{equation}
 E_L\leq e_0+E_{L-1} \Rightarrow E_L\leq e_0+e_1+E_{L-2}\leq \sum_{j=0}^{L-1}e_j+E_0=\sum_{j=0}^{L-1} e_j.
\end{equation}

The local truncation error $e_j$ has two parts -- one shown in \cref{thm:lte-mag2} and the other from the numerical quadrature. Note that the numerical quadrature of Riemann sum follows the standard numerical analysis, see e.g., \cite[Theorem 2]{FangLiuSarkar2024} and \cite{BurdenNA}. Specifically, we have
\begin{equation}
    \norm{U_2(t_{j+1}, t_j)  - \Tilde{
U}_2(t_{j+1}, t_j)} \leq \frac{\Delta t^2}{M}\sup_{s\in [t_j,t_{j+1}]} \norm{H'(s)}+\frac{3 \Delta t^3}{M} \sup_{s\in [t_j,t_{j+1}]} \norm{H(s)} \norm{H'(s)} .
\end{equation}
Combining both parts, we obtain the desired result.
\end{proof}

We make two remarks here. First, we can further bound $\beta_\mathrm{comm}^j$ by $\alpha_\mathrm{comm}^j$. This is because
\begin{equation} 
\begin{aligned}
 & \sup_{\sigma, s, t, \tau\in [t_j,t_{j+1}]}\|  \left[[H(\sigma),H(s)],[H(\tau), H(t)]\right] \|
  \\
  \leq & 4 \sup_{\tau\in [t_j,t_{j+1}]} \norm{H(\tau)} \sup_{\sigma, s, t\in [t_j,t_{j+1}]} \|\left[[H(\sigma),H(s)],H(t)\right]\|.
\end{aligned}
\end{equation}
Therefore, the long-time bound can be rewritten as
\begin{equation}\label{eq:long_time_with_quadrature}
\begin{aligned}
    & \norm{U_\mathrm{exact}( T, 0) - \tilde{U}_2(T, 0)} 
    \\
    \leq & \Delta t^3 \sum_{j=0}^{L-1}\alpha_\mathrm{comm}^j 
    \left( \frac{13}{24} + \frac{5}{12}  \Delta t \alpha \right) +\frac{T \Delta t}{M} (1+3 \Delta t \alpha)\sup_{s\in [0,T]}\norm{H'(s)}.
    \end{aligned}
\end{equation}
Second, it is important to note that while \cref{eq:long_time_with_quadrature} contains a term that is linear in $\Delta t$, this does not imply the algorithm is only of first order. This is due to the freedom to choose $M$ in the quantum algorithm and the low cost of implementing the quadrature via quantum circuits. In particular, $M$ is chosen such that the last term in \cref{eq:long_time_with_quadrature} is of comparable magnitude to the first two terms. In other words, we choose 
\begin{equation}    \label{eq:M_estimate}
M = \Or\left(\sup_{s\in[0,T]}\norm{H'(s)} \frac{T(1+3\Delta t \alpha)}{\Delta t^2 \sum_{j=0}^{L-1}\alpha_\mathrm{comm}^j (13+10 \Delta t \alpha)}\right)
\end{equation}
in the quantum algorithm, and the circuit gate complexity due to LCU is only $\Or(\log(M))$, which contributes to a logarithmic dependence on the derivative of the Hamiltonian and the time $T$ as discussed in \cref{sec:revisit_quantum_algorithm}.

\section{Superconvergence in the Discrete Case} \label{sect:superconv}
In this section, we discuss in detail how to establish the superconvergence with a finite $N$ rigorously for the unbounded Hamiltonian simulation task in the interaction picture. We consider the interaction picture Hamiltonian~\cref{eq:Ham_interaction_picture} with $A$ and $B$ defined as \cref{eq:A_fd} and \cref{eq:B_fd} with a finite grid number $N$.
Thanks to \cref{thm:long-time-mag2}, it is sufficient to show that the nested commutator
\begin{align}
   & \sup_{\sigma, \tau, t \in [0, \Delta t]} \norm{ [H(\sigma),[H(\tau), H(t)]]} 
   \\
   = &  \sup_{\sigma, \tau, t \in [0, \Delta t]} \norm{ \left[ e^{iA\sigma} B e^{-iA\sigma},  [e^{iA\tau} B e^{-iA\tau}, e^{iAt} Be^{-iAt} ]\right]} \leq C \Delta t^2,
\end{align}
with a constant $C$ independent of $N$. Note that for any unitary $U$, the commutator has the property
\begin{equation}
    U^\dagger [O_1, O_2] U =  [U^\dagger O_1 U, U^\dagger O_2 U],
\end{equation}
\begin{equation}
    U^\dagger [O_1, [O_2, O_3]] U =  [U^\dagger O_1 U, [U^\dagger O_2 U, U^\dagger O_3 U]] ,
\end{equation}
for any operators $O_1, O_2, O_3$. By choosing $U = e^{-iAt}$, we can see that it is sufficient to estimate
\begin{equation}
   \sup_{\sigma, \tau, t \in [0, \Delta t]}   \norm{ \left[ e^{iA(\sigma-t)} B e^{-iA(\sigma-t)}, [e^{iA(\tau-t)}B e^{-iA(\tau-t)}, B] \right]}\leq C \Delta t^2,
\end{equation}
with a constant $C$ independent of $N$.
In other words, we need to establish the following theorem.

\begin{thm} \label{thm:comm_estimate_dis}
Let $V\in S(1)$, and $A$ and $B$ are the finite difference discretizations of $-\frac{1}{2}\Delta$ and $V(x)$ with periodic boundary conditions via $N$ grid points.
Then for $1/N \le \Delta t \le1$, 
\begin{equation*}\label{eq:mainineq}
\sup_{s, \tau \in [-\Delta t, \Delta t]} \left\|\left[ \left[B,e^{i\tau A} B e^{-i \tau A}\right],  e^{is A} B e^{-is A}
 \right]\right\|\le C \Delta t^2
\end{equation*}
with the constant $C$ independent of $N$ and $\Delta t$. Here $\norm{\cdot}$ denotes the matrix operator norm.
\end{thm}

\REV{Here $V \in S(1)$, meaning that $V$ is bounded together with all its derivatives. The definition of $S(1)$ is given in \cref{eq:def_S(1)}. We also note that this theorem seems to require an additional condition on the step size $\Delta t$ and $N$, namely $1/N \leq \Delta t$. However, we emphasize that this is not an additional restriction. Allowing a larger time step is precisely what we want here. We would like to show that under the superconvergence property, we can still take a time step $\Delta t$ that is large (and not restricted by the size of $N$), so that the number of time steps required in the algorithm can be chosen independently of $N$. In this way, the circuit cost per step is $\polylog(N)$, and we can maintain the overall algorithmic cost at $\polylog(N)$. In other words, as $N$ increases, the condition does not require $\Delta t$ to decrease accordingly for the bound to hold, unlike the case for the Taylor argument. Note that this corresponds to the opposite regime from that in which the Taylor series argument would work. For the Taylor expansion to hold, one needs to require the time step $\Delta t$ to be very small -- particularly, in our case, much smaller than $1/N$. This would require the total number of time steps $L = T/\Delta t$ to be $\poly(N)$, which would fail to yield the overall $\polylog(N)$ cost that we aim for.}

The rest of this section is dedicated to rigorously establishing this theorem. We begin by reformulating the discrete matrix spectral norm estimation problem as estimating the operator norm of a continuous (in space) operator in the  $L^2$ to $L^2$ sense, using discrete microlocal analysis techniques previously introduced to analyze Trotterization in \cite{BornsWeilFang2022}.
However, this alone does not fully resolve the issue due to a scaling mismatch that prevents the direct application of Egorov's theorem. To address this, we apply a unitary transformation to the problem and leverage a more advanced framework—the two-parameter symbol class~\cite{NonnenmacherSjostrandZworski2014}—to overcome this difficulty. To the best of our knowledge, this is the first application of the two-parameter symbol class in both numerical analysis and quantum algorithms.

\subsection{Converting from Matrix to Operator and Challenges}\label{sec:finite_difference_discretization}
In our previous work~\cite{BornsWeilFang2022}, an analysis of the Trotterization of the semiclassical Schrödinger equation is provided.
While the topic is entirely different from this work, we introduced discrete microlocal analysis~\cite{ChristiansenZworski2010,Schenck2009,Deleporte2019,DyatlovJezequel2021} -- the Weyl quantization on a torus -- as a tool for numerical analysis.

The key is that such quantization for a phase-space function $f$ corresponds to an operator acting on the $N$-dimensional linear space, i.e. an $N$ by $N$ matrix. The Weyl quantization is denoted as $\op_N(f (x,\xi))$, and the $N$-dimensional linear space is denoted as $H_N$.
There are two important results established in \cite[Section 4.2]{BornsWeilFang2022} regarding the discrete microlocal analysis that are crucial for the numerical analysis purpose here.
First, the matrix that appears in the finite-difference discretization of $-\frac{1}{2}\Delta$, namely,
\begin{equation}
\tilde{A}_0: = \frac{1}{2}\left(\begin{array}{ccccc}
        2 & -1 & & & -1 \\
         -1& 2 & -1 & & \\
          & \ddots& \ddots& \ddots& \\
           & & -1& 2 & -1\\
        -1& & & -1 & 2 \\
    \end{array}\right)_{N \times N} = \op_N( 1- \cos(\xi)).
\end{equation}
is the Weyl quantization of a $S(1)$ symbol on the torus. For $d > 1$, we analogously obtain $\text{op}_N \left( \sum_{j=1}^d (1 - \cos\xi_j) \right)$. Similarly, $B$, as defined in \cref{eq:B_fd}, equals $\op_N(V(x))$.
Second, it has been shown that when taking $h = \frac{1}{2 \pi N}$, the spectral norm of the pseudodifferential operators on $H_N$ are upper bounded by the norm from $L^2(\mathbb{R}^d)$ to $L^2(\mathbb{R}^d)$ of the (standard) continuous Weyl quantization of the phase-space function. This result can be summarized as follows.

\begin{lemma}[Lemma 9 in \cite{BornsWeilFang2022}]\label{lem:H_N_L^2_bound}
Let $a\in\cinf(\TT^{2d})$, and let $\tilde a$ be the lift of $a$ to a periodic function on $T^*\RR^d$. Then
\begin{equation}\label{eq:H_N_norm_bound}
\|\op_N(a)\|_{H_N\to H_N}\le\|\op_h(\tilde{a})\|_{L^2(\RR^d)\to L^2(\RR^d)}.
\end{equation}
\end{lemma}

In light of this, we can work with the corresponding (continuous) operators with their spatial degrees of freedom kept continuous. In particular,
note that
\begin{equation}
    A = \frac{N^2}{(b-a)^2} \tilde{A}_0 = \frac{1}{h^2}\frac{1}{(2\pi)^2(b-a)^2} \op_N( \tilde{a}_0),
\end{equation}
where $\tilde{a}_0 \in S(1) $. 
Thus, to understand $A$, we can consider the operator
\begin{equation}
    \frac{1}{h^2} A_0, \quad A_0 = \op_h(a_0), \quad a_0 \in S(1),
\end{equation}
where we move the absolute constants into the definition of $a_0$. Also, $B = \op_N(V(x))$ corresponds to consider 
\begin{equation}
    \op_h(V(x)) = V(x).
\end{equation}
Therefore, to obtain our main result~\cref{thm:comm_estimate_dis}, we need to estimate
\begin{equation}\label{eq:main_comm_estimate_cts}
\sup_{\tau, s\in [-\Delta t, \Delta t]}
\left\lVert 
\left[ \left[ V(x), e^{i \tau \frac{A_0}{h^2}} V(x) e^{-i \tau \frac{A_0}{h^2}}\right] ,  e^{i s \frac{A_0}{h^2}} V(x) e^{-i s \frac{A_0}{h^2}}  \right] \right\rVert_{L^2(\mathbb{R}^d) \to L^2(\mathbb{R}^d)}.
\end{equation}

Note that in previous continuous superconvergence proof, the crucial step is the exact Egorov's theorem, namely,
\begin{equation}
e^{ - i \frac{\tau}{2}\Delta  } V(x) e^{ i \frac{\tau}{2}\Delta} = \op_{h_0}(V(x + p\tau)),
\end{equation}
with $h_0 = 1$.
This holds exactly because here the unitaries are governed by a quadratic symbol $\xi^2$, which is no longer the case for our scenario -- the unitaries in \cref{eq:main_comm_estimate_cts} are governed by an $S(1)$ symbol instead. While the general version of Egorov's theorem still holds, it is only valid up to the Ehrenfest time $\Or(\log(h^{-1}))$. However, in our setting, it appears that we need to consider a much longer time scale $\Or(h^{-2})$. This is the main challenge. We emphasize that the discrete superconvergence result is not merely a direct consequence of discrete microlocal analysis -- new techniques are needed.

\subsection{Result of Continuous Symbol Calculus}

\cref{sect:symbol_calculus} will be focused on proving the following theorem. In this section, we discuss the algorithmic implications of this result. \REV{We refer the readers to the definition of $S(1)$ in \cref{eq:def_S(1)}.}

\begin{thm}\label{thm:main_S(1)symbol}
Let $v \in S(1)$, and let $V$ be multiplication by $v(x)$. Let $A_0 = \text{op}_h(a)$ for $a\in S(1)$ depending only on $\xi$.
Then, for $0 \leq \Delta t \leq 1$ and $h \leq \Delta t \leq 1$,

\[
\sup_{\tau, s\in [-\Delta t, \Delta t]}
\left\lVert \left[ \left[ V, e^{i \tau \frac{A_0}{h^2}} V e^{-i \tau \frac{A_0}{h^2}}\right] ,  e^{i s \frac{A_0}{h^2}} V e^{-i s \frac{A_0}{h^2}}  \right] \right\rVert_{L^2(\mathbb{R}^d) \to L^2(\mathbb{R}^d)} \leq C_V \Delta t^2,
\]
where $C_V$ is a constant depending only on $V(x)$ and the dimension $d$, uniformly in $h$ and $\Delta t$.
\end{thm}

An immediate consequence of \cref{thm:main_S(1)symbol} is \cref{thm:comm_estimate_dis}, thanks to \cref{lem:H_N_L^2_bound}. Thus, the nested commutators as defined in \cref{eq:nest_comm_alpha} and \cref{eq:nest_comm_beta} can be estimated as
\begin{equation}
    \alpha_\mathrm{comm}^j  =  C_V \Delta t ^2
 , \quad  \beta_\mathrm{comm}^j  = 4 C_V  \sup_{\tau \in [t_j,t_{j+1}]} \norm{H_I(\tau)}\Delta t ^2,
\end{equation}
where $C_V$ and $C'_V$ are constants depending on $V$ but independent of $N$.
Combining with \cref{eq:long_time_with_quadrature}, we have 
\begin{equation}\label{eq:long_time_with_quadrature_superconv}
\begin{aligned}
    \norm{U_\mathrm{exact}( T, 0) - \tilde{U}_2(T, 0)} 
    &\leq \Delta t^4 T C_V 
    \left( \frac{13}{24} + \frac{5}{12}  \Delta t \alpha \right) +C \frac{T \Delta t}{M} (1+3 \Delta t \alpha) N,
    \end{aligned}
\end{equation}
with an absolute constant $C$, 
where we used the fact that $\norm{H'(t)} = \norm{[A,B]} = \Or(N)$.
This completes the proof of fourth-order superconvergence, with the error preconstant independent of $N$. As discussed in \cref{eq:M_estimate}, we choose $M$ so that the second term is comparable in size to the first term, while the quantum gate complexity remains only $\log(M)$.

\REV{\subsection{Quantum Algorithm Complexity}}
\REV{
In this section, we discuss the error analysis and the corresponding implications for the quantum algorithm’s complexity.
As shown in \cite[Theorem~12]{FangLiuSarkar2024}, the relationship is given by:
\begin{lemma}[Long-time complexity of the second-order quantum Magnus algorithm]
\label{lem:long_time_cost}
Let the Hamiltonian $H(s)$ be accessed via a block-encoding denoted as HAM-T with a subnormalization factor $\alpha$ and $n_a$ ancillas (see \cite[Section 5.1]{FangLiuSarkar2024} for the concrete definition).
Assume that
\begin{equation}
    \|U_{\mathrm{exact}}(t_{j+1}, t_j) - U_2(t_{j+1}, t_j)\| \leq C_H \Delta t^{1+\theta},
\end{equation}
for some non-negative real number $\theta$ and a constant $C_H$ that may depend on $H$, and $\Delta t = T/L$ is the short time step with the total number of time steps denoted as $L$.
Then, for any precision $0 < \delta < 1$ and $T > \delta$, the quantum Magnus algorithm can implement an operation $W$ such that 
\[
\|W - U_\mathrm{exact}(T,0)\| \leq \delta
\]
with failure probability $\mathcal{O}(\delta)$, and with the following computational cost:
\begin{enumerate}
    \item $\mathcal{O}\!\left(\alpha T + \frac{C_H^{1/\theta} T^{1+1/\theta}}{\delta^{1/\theta}} 
    \log\!\left(\frac{C_H T}{\delta}\right)\right)$ uses of $\text{HAM-T}$, its inverse, or its controlled version;
        \item $\mathcal{O}\!\left(\! (n_a + \log(M)) \left(\alpha T + 
    \frac{C_H^{1/\theta} T^{1+1/\theta}}{\delta^{1/\theta}} 
    \log\! \left( \frac{C_H T}{\delta} \right) \right)\!\right)$ elementary gates, where $M$ is the number of numerical quadrature points used in the algorithm;
    \item $\mathcal{O}\!\left(\log\!\frac{\max_{s\in [0,T]} \|H'(s)\| T}
    {C_H \delta}\right)$ ancilla qubits.
\end{enumerate}
\end{lemma}
}

\REV{
We now consider the algorithmic complexity for our case when superconvergence is achieved, in particular to justify that the overall cost scales polylogarithmically with respect to $N$, rather than polynomially in $N$.
In our setting, $\text{HAM-T}$ corresponds to the interaction-picture Hamiltonian (\cref{eq:Ham_interaction_picture}), which can be constructed following \cite[Fig.~3]{FangLiuSarkar2024}. In our case, $e^{-isA}$ can be diagonalized by the quantum Fourier transform (QFT), and hence can be implemented via a QFT, controlled rotations, and an inverse QFT (IQFT). The overall gate complexity for constructing $\text{HAM-T}$ is $\mathcal{O}(\polylog (N))$.
The subnormalization factor $\alpha$ in our case is $\norm{B} = \norm{V}_\infty$, which is a bounded constant (not depending on $N, T, \delta$), where $V(x)$ denotes the potential. Since $B$ is diagonal, the block-encoding requires only one ancilla qubit. 
}

\REV{
Moreover, the local truncation error in our case is given by \cref{thm:lte-mag2}. Combining this with our main superconvergence estimate \cref{thm:comm_estimate_dis}, we have
$C_H = C_V$, where $C_V$ is a constant defined in \cref{thm:comm_estimate_dis} and does not depend on $N$, $\Delta t$, $T$ or $\delta$. We also have $\theta = 4$, which is the convergence order of the algorithm.
}

\REV{
To achieve an overall error tolerance of $\delta$ in the full evolution, we choose the short time step $\Delta t$ by balancing the quadrature error and the Magnus error (see \cref{eq:long_time_with_quadrature} and the discussion that follows).
The idea is to choose $M$ such that the quadrature error and the Magnus error are comparable, and then choose $\Delta t$ (or equivalently, the number of time steps $L$) such that the Magnus error is of order $\delta$ (up to a constant, say $1/2$).
The total error is then controlled by $\Delta t$ via the triangle inequality (see \cite[Section~5]{FangLiuSarkar2024} for a detailed discussion).
In light of this, we determine the total number of time steps $L$ by choosing the time-step size $\Delta t = T/L$ such that
\begin{equation}
    \Delta t \leq \left(\delta /(C_H T)\right)^{1/4},
\end{equation}
so we can choose $\Delta t = \Theta \left(\delta^{1/4} / T^{1/4}\right)$.
We then set $M$ to be
\begin{equation}
    M = \Or (N /\Delta t^3) = \Or(N T^{3/4}/ \delta^{3/4}),
\end{equation}
where we have used the fact that $\norm{H'(s)} = \norm{e^{iAs} [A,B] e^{-iAs}} = \norm{ [A,B]} = \Or(N)$.
}

\REV{
Therefore, the query complexity of HAM-T is 
\begin{equation}
    \Or\left (\frac{T^{5/4}}{\delta^{1/4} } \log\left(\frac{T}{\delta}\right) \right).
\end{equation}
and the elementary gates needed to construct the algorithm, other than applications of HAM-T are
\begin{equation}
    \Or \left(  \log(N) \frac{T^{5/4}}{\delta^{1/4} } \log^2\left(\frac{T}{\delta}\right) \right).
\end{equation}
Including the number of elementary gates required to construct the $\text{HAM-T}$ block-encoding, the overall gate complexity is estimated as
\begin{equation}
    \Or\left (\polylog(N) \frac{T^{5/4}}{\delta^{1/4} } \log^2\left(\frac{T}{\delta}\right) \right).
\end{equation}
Note that we do not aim to optimize the $\polylog(N)$ factor here, but it is primarily determined by the number of elementary gates required for the QFT, which, according to \cite{HalesHallgren2000}, is $\mathcal{O}(\log N \log\log N)$.
}

\section{Symbol Calculus and Proof of \cref{thm:main_S(1)symbol}} \label{sect:symbol_calculus}
In this section, we prove \cref{thm:main_S(1)symbol}. The key idea is to apply a unitary transformation and identify the resulting operators as the Weyl quantization of two-parameter symbols with appropriately defined semiclassical parameters.

\subsection{Standard and Two-parameter Symbol Calculus}\label{sec:subsection_revisit_symbol_class}
In this section, we revisit only necessary concepts from semiclassical analysis. For a more detailed introduction, please see textbooks such as \cite{Martinez2002,zworski2022semiclassical}. Let $\epsilon \in (0, 1]$. For a phase-space function $a(x,\xi)$ in a symbol class, we can define its Weyl quantization $\op(a) = \op_\epsilon(a): \mathcal{S}(\mathbb{R}^d) \to \mathcal{S}(\mathbb{R}^d)$ by
\begin{equation}\label{eq:weylquant}
(\op_\epsilon(a)u)(x):=\frac{1}{(2\pi \epsilon)^d}\iint_{\RR^{2d}} e^{\frac{i}{\epsilon}\la x-y,\xi\ra}a\left(\frac{x+y}{2},\xi\right)u(y)\,dy\,d\xi.
\end{equation} 
The Weyl quantization is well-defined for a variety of symbol classes. For us, the useful one is the $S(1)$ symbol class, defined as follows: \begin{equation} \label{eq:def_S(1)}
    S(1)=S(1; \RR^{2d}):=\left\{a\in\cinf(T^*\RR^d):\left|\dd_x^{\alpha}\dd_{\xi}^{\beta}a(x,\xi)\right|\le C_{\alpha\beta}\right\}.
\end{equation}
For $S(1)$ symbols, the Weyl quantization $\op_\epsilon(a)$ is in fact bounded from $L^2(\mathbb{R}^d)$ to $L^2(\mathbb{R}^d)$, thanks to the Calderon-Vaillancourt theorem.

For the $S(1)$ symbols, one also has the Egorov's theorem (see \cite[Theorem 11.1 and Remark (ii)]{zworski2022semiclassical} for the proof).
Suppose $U = U(1)$, where $U(t)$ is the propagator of 
\begin{equation}\label{eq:schd_S(1)_Q}
    i \epsilon \partial_tU(t) = Q U(t),\quad Q=\Op_\epsilon(q),
\end{equation}
with $U(0) = I$ and $q \in S(1)$. Then we have, for $A = \op_\epsilon(a)$ with $a \in S(1)$ that
\begin{equation}\label{eq:egorov_S1}
    U^{-1} A U = B, \quad B = \op_\epsilon(b),\quad b - \kappa^\ast a \in \epsilon S(1), 
\end{equation}
where $\kappa^\ast = \kappa^\ast_1$ is the symplectomorphisms defined by $q_t$ at the time $t = 1$.

Another important symbol class for the discrete superconvergence proof is the two-parameter symbol class. Here we follow the notations as used in \cite{NonnenmacherSjostrandZworski2014}. Specifically, the \emph{two-parameter symbol class} with respect to the semiclassical parameters $\epsilon$ and $\tilde{\epsilon}$ are defined as
\begin{equation} \label{S_1/2}
a \in  \tilde{S}_{1/2} : = \tilde{S}_{1/2}(1; T^*\mathbb{R}^d) \iff |\partial^\alpha a| \leq C_\alpha \left(\tilde{\epsilon}/\epsilon \right)^{\frac{1}{2}|\alpha|}, \quad \epsilon, \tilde{\epsilon} \in [0,1].
\end{equation}
\REV{Here the subscript $\frac{1}{2}$ refers to the prefactor in front of $|\alpha|$ in the power of the right hand side. It is important to note that unlike the standard symbol class $S(1)$ defined in \cref{eq:def_S(1)},  $\tilde{S}_{1/2}$ allows the symbol to have certain dependence on two semiclassical parameters $\epsilon$ and $\tilde{\epsilon}$. For notational simplicity, we adopt the convention as in \cite{NonnenmacherSjostrandZworski2014} and suppress the dependence of $a$ on $\epsilon$ and $\tilde \epsilon$ in its notation. The (main) semiclassical parameter $\epsilon$ is still the parameter used to define the Weyl quantization as in \cref{eq:weylquant}, while $\tilde \epsilon$ is an additional scale that we will utilize to tune later in our proof. The role of the two-parameter symbol class will be explained in \cref{sec:necessity_two_parameter_symbol}.}
The corresponding classes of operators by $\Psi_{1/2}^{m, \tilde{m}, k}(\mathbb{R}^d)$ or $\tilde{\Psi}_{1/2}$. 
Correspondingly, we also have the Weyl quantization defined in the same manner for this two-parameter symbol class. 

One of the fundamental properties of the two-scale symbol class is the composition rule (\cite[Lemma 3.2]{NonnenmacherSjostrandZworski2014}), which follows a similar form in the $ \tilde{S}_{1/2} $ class as in the $ S(1) $ case.
Specifically, if $a, b \in \tilde{S}_{1/2}$ and $\op_\epsilon(c) = \op_\epsilon(a) \circ \op_\epsilon(b)$, then for any integer $N > 0$, $c$ reads
\begin{equation}\label{eq:composition1_2}
c(x, \xi) = \sum_{k=0}^{N} \frac{1}{k!} \left( \frac{i \epsilon}{2} \sigma(D_x, D_\xi; D_y, D_\eta) \right)^k a(x, \xi) b(y, \eta) \Big|_{x=y,\xi=\eta} + e_N(x, \xi),
\end{equation}  %
where the remainder $e_n$ is bounded by $\epsilon^{N+1}$. 
Similarly, the Calder\'on-Vaillancourt Theorem naturally extends to the symbol class $\tilde{S}_{1/2}$ so that the Weyl quantization of a $\tilde{S}_{1/2}$ symbol is bounded from $L^2$ to $L^2$. For a detailed discussion and the proof of such properties, see, e.g., \cite[Section 3]{NonnenmacherSjostrandZworski2014} and \cite[Chapter 6]{Nonnenmacher_book}.
Additionally, for $\tilde{S}_{1/2}$, we also have the corresponding Egorov's theorem as in \cite[Lemma 3.12]{NonnenmacherSjostrandZworski2014}:
for $A = \op_{\epsilon}(a)$, $a \in \tilde{S}_{1/2}$, we have
\begin{equation} \label{Egorov S_{1/2}}
U^{-1} A U = B, \quad B = \op_{\epsilon}(b), \quad b - \kappa^* a \in {\epsilon}^{\frac{1}{2}} \tilde{\epsilon}^{\frac{3}{2}} \tilde{S}_{1/2}.
\end{equation}

\REV{\subsection{The Necessity of Two-parameter Symbol Class and Proof Intuition} \label{sec:necessity_two_parameter_symbol}}

\REV{
In the previous section, we introduced both the standard and two-parameter symbol classes. The reason for introducing the latter may not be immediately clear. In this section, we explain its necessity, together with the intuition behind the proof.}

\REV{
For illustration purposes, we consider one layer of the commutator instead of the two layers of nested commutators in \cref{thm:main_S(1)symbol} for simplicity, namely,
\begin{equation}\label{eq:one-layer-comm}
    \left[V, e^{i  \frac{s}{h^2} A_0} V e^{-i \frac{s}{h^2}A_0} \right].
\end{equation}
Intuitively, we hope that every layer of nested commutator contributes another factor of $s$, which would lead to $s^2$ and taking the supremum over $s \in [-\Delta t, \Delta t]$ yields $\Delta t^2$ as desired for two layers, as expected in \cref{thm:main_S(1)symbol}.
}

\REV{
In our original numerical analysis problem, there are two relevant scales: one is the time scale $s$ (for the short time step), and the other is the spatial discretization scale $1/N = h$. Hence, we would expect two parameters to matter. More concretely, one can observe that in~\cref{eq:one-layer-comm} there are two scales denoted by slightly different parameters that actually matter: one is the time scale $s/h^2$ which resembles the reciprocal of the semiclassical parameter $\epsilon$ in the Egorov theorem as in \cref{eq:egorov_S1}, but in order to apply this standard Egorov theorem with respect to the standard symbol class, one needs the reciprocal of the time scale in unitary evolution and the semiclassical parameter in the Weyl quantization to match. Namely,
the unitary in the Egorov theorem needs to be $e^{-i\frac{1}{\epsilon}\op_\epsilon (q)} $, but here we have a mismatch, where the semiclassical parameter in the Weyl quantization is $h$, while the reciprocal of the time scale is $h^2/s$. This observation motivates the definition of a new characteristic parameter $h' := h^2/s$.
}

\REV{
To apply the Egorov theorem (of any kind -- either the standard one-parameter or more exotic two-parameter), we need to make these parameters coincide. The question is how. It turns out that one can define a unitary operator $U_\alpha$ on $L^2(\mathbb{R}^d)$ by 
\begin{equation}
    U_\alpha f(x): = \alpha^{d/2} f(\alpha x).
\end{equation} 
This is indeed a unitary transformation, as proved in \cref{lem:def_u_unitary} (for some real-valued $\alpha$ to be determined). By choosing $\alpha = \frac{h}{s}$, we have $U_\alpha^\dagger \op_h(a) U_\alpha = \op_{h'} (a)$ (see \cref{lem:UAU}).
This is very helpful, as we can multiply $U_\alpha$ and $U_\alpha^\dagger$ to \cref{eq:one-layer-comm} without changing its operator norm. We get
\begin{equation}\label{eq:comm_after_U_alpha}
 U_\alpha^\dagger \left[V, e^{i  \frac{s}{h^2} A_0} V e^{-i \frac{s}{h^2}A_0} \right] U_\alpha
 =  
 \left[U_\alpha^\dagger V U_\alpha, e^{i  \frac{s}{h^2} U_\alpha^\dagger A_0 U_\alpha}  U_\alpha^\dagger V U_\alpha e^{-i \frac{s}{h^2} U_\alpha^\dagger A_0 U_\alpha} \right]. 
\end{equation}
In order to match the time evolution with the correct scaling as explained above, we need to have 
\begin{equation}\label{eq:UA0U}
    U_\alpha^\dagger A_0 U_\alpha =  \op_{h'} (a), \quad h' = h^2/s.
\end{equation}
This uniquely determines the choice of $\alpha$ and its corresponding unitary $U : = U_\alpha$, with $\alpha = h/s$. This is the unitary that we consider in~\cref{sec:subsection_unitary}.
}

\REV{
But under this unitary transformation, we can see that the commutator in \cref{eq:comm_after_U_alpha} now concerns the operator 
$U_\alpha^\dagger V U_\alpha = U^\dagger  V U$, which we prove in \cref{lem:UVU} that it can be identified with the symbol $v(\frac{s}{h} x)$. But importantly, it is \textbf{no longer a $S(1)$ symbol}. One can check its derivatives and observe that they are not uniformly bounded in terms of $s$ and $h$ but rather depend on such parameters. So it is no longer the standard one-parameter symbol class that one typically works with (see \cref{eq:def_S(1)}). The corresponding standard one-parameter symbol calculus and Weyl quantization no longer apply to this problem of our interest. However, the two-parameter symbol class comes to the rescue. Notice that we can still show that 
$v(\frac{s}{h}x) \in  \tilde{S}_\frac{1}{2}$ with respect to the parameter $h' = h^2/s$ and $s$. In particular, observe that every partial derivative acting on $v$ leads to a scaling of 
\begin{equation}
    \frac{s}{h} = \sqrt{\frac{s^2}{h^2}} = \sqrt{\frac{s}{h'}},
\end{equation}
where in the two-parameter symbol definition~\cref{S_1/2}, we can choose $\epsilon = h' = h^2/s$ and $\tilde \epsilon = s$. This alignment is ideal, as the main semiclassical parameter $\epsilon = h'$ coincides with the semiclassical parameter in the unitary evolution for \cref{eq:UA0U}. Without this alignment, even the two-parameter symbol class would not effectively support the argument. It now becomes evident that we need the Weyl calculus for such a two-parameter symbol class, including the commutator rules,  the corresponding two-parameter Egorov-type theorems, and the Calderon-Vaillancourt theorem, which luckily were discussed and established in detail in \cite{NonnenmacherSjostrandZworski2014}. For self-containedness, we revisit such results in \cref{sec:subsection_revisit_symbol_class}. It is worth noting that the leading terms in the two-parameter commutator rules (\cref{eq:composition1_2}) include powers of the (main) semiclassical parameter $\epsilon$, and the Egorov theorem (\cref{Egorov S_{1/2}}) involves both. It is still crucial to show that after a careful application of this more exotic symbol calculus, all scalings are well-behaved and uniformly bounded by $\Delta t^2$ as we desired, which is carefully carried out in \cref{sec:subsection_proof_main_thm}. It is also nice to note that $\epsilon  = h' = h^2/s \leq s$, if we have $h \leq s$, as in the condition of \cref{thm:main_S(1)symbol}, and also that the other semiclassical parameter $\tilde{\epsilon} = s$, so both semiclassical parameters under the condition of \cref{thm:main_S(1)symbol} can be bounded by $s$ uniformly (in $h$). The discussion of this section is for illustration purposes and we only consider a single layer of the commutator, in two-layer case, we need to set $\Delta t : = h_0$ and rescale the corresponding two time-scales properly, see the proof in \cref{sec:subsection_proof_main_thm}.
}

\subsection{Unitary Transformation and Properties}\label{sec:subsection_unitary}

\begin{lemma}\label{lem:def_u_unitary}
Let $U$ be an operator on $L^2(\mathbb{R}^d)$ defined as
\begin{equation}\label{eq:u_def}
    Uf(x)=\left(\frac{h}{s} \right)^\frac{d}{2} f\left(\frac{h}{s}x\right).
\end{equation} 
Then $U$ is a unitary operator.
\end{lemma}
\begin{proof}
    It suffices to show that $U^\dagger U  = U U^\dagger = I$. 
Given definition of the operator $ U $ acting on $ f $, we can use the fact that $ U^\dagger $ is the adjoint of $ U $ to find $ U^\dagger $, which satisfies
\begin{equation}
\langle Uf, g \rangle_{\REV{L^2}} = \langle f, U^\dagger g \rangle_{\REV{L^2}} \quad \text{for all} \quad f, g \in L^2(\mathbb{R}^d).
\end{equation}
First, we compute the inner product $ \langle Uf, g \rangle_{\REV{L^2}} $. Using the definition of $ U $
\begin{equation}
\langle Uf, g \rangle_{\REV{L^2}} = \int_{\mathbb{R}^d} U f(x) \overline{g(x)} \, dx = \int_{\mathbb{R}^d} \left( \frac{h}{s} \right)^{d/2} f\left( \frac{h}{s} x \right) \overline{g(x)} \, dx.
\end{equation}
We use the change of variables $ y = \frac{h}{s} x $, so $ x = \frac{s}{h} y $ and $ dx = \left( \frac{s}{h} \right)^d dy $. Since the scaling factor for the volume element in $ d $-dim is $ \frac{s}{h} $, so on $ \mathbb{R}^d $, the total scaling factor for the volume element in $ d $-dim is the product of scaling factors in each dim, which gives $ \left( \frac{s}{h} \right)^d $. This product corresponds to the determinant of the Jacobian matrix for the transformation. \REV{It can be easily check that $U^\dagger U = U U^\dagger = I$, for this construction of $U$. For completeness, we present this in \cref{pf:appendix_UdagU}.}
\end{proof}

\begin{lemma} \label{lem:UAU}
Suppose that $a(\xi) \in S(1)$ is a function of only $\xi \in \mathbb{R}^d$. For $A_0 = \op_h (a)$  and $U$ as defined in \cref{eq:u_def}, we have
    \begin{equation}
        U^\dagger A_0U
        =\op_{\frac{h^2}{s}}a.
    \end{equation}
\end{lemma}
\begin{proof}
    First, we apply $ U $ to a function $ f(x) $
    \begin{equation}
    (Uf)(x) = \left( \frac{h}{s} \right)^{d/2} f\left( \frac{h}{s} x \right).
    \end{equation}
    We apply $ A_0 := \op_h(a) $ to $ Uf(x) $, and using the definition of Weyl quantization, we have
    \begin{equation}
    (\op_h(a) (Uf))(x) = \frac{1}{(2\pi h)^d} \iint_{\mathbb{R}^{2d}} e^{i \langle x - y, \xi \rangle / h} a(\xi ) (Uf)(y) \, dy \, d\xi.
    \end{equation}
    Substituting $ Uf(y) = \left( \frac{h}{s} \right)^{d/2} f\left( \frac{h}{s} y \right) $ give us
    \begin{equation}
    (\op_h(a) (Uf))(x) = \left( \frac{h}{s} \right)^{d/2} \frac{1}{(2\pi h)^d} \int_{\mathbb{R}^{d}}\int_{\mathbb{R}^{d}} e^{i \langle x - y, \xi \rangle / h} a(\xi) f\left( \frac{h}{s} y \right) \, dy \, d\xi.
    \end{equation}
    We then use the change of variables: Let $ z = \frac{h}{s} y $, then $ dy = \frac{s^d}{h^d} dz $, and the expression becomes
    \begin{equation}
    (\op_h(a) (Uf))(x) = \left( \frac{h}{s} \right)^{d/2} \frac{1}{(2\pi h)^d} \int_{\mathbb{R}^{d}}\int_{\mathbb{R}^{d}} e^{i \langle x - \frac{s}{h} z, \xi \rangle / h} a(\xi) f(z) \frac{s^d}{h^d} \, dz \, d\xi.
    \end{equation}
    Simplifying it yields
    \begin{equation}
    (\op_h(a) (Uf))(x) = \left( \frac{s}{h} \right)^{d/2} \frac{1}{(2\pi h)^d} \int_{\mathbb{R}^{d}}\int_{\mathbb{R}^{d}} e^{i \langle x - \frac{s}{h} z, \xi \rangle / h} a(\xi) f(z) \, dz \, d\xi.
    \end{equation}
    Now we apply $ U^\dagger $, which is defined as
    \begin{equation}
    U^\dagger g(x) = \left( \frac{s}{h} \right)^{d/2} g\left( \frac{s}{h} x \right).
    \end{equation}
    We apply this to $ (\op_h(a) (Uf))(x) $ to get
    \begin{equation}
    (U^\dagger (\op_h(a) (Uf)))(x) = \left( \frac{s}{h} \right)^{d/2} (\op_h(a) Uf)(x).
    \end{equation}
    Substituting the expression for $ (\op_h(a) (Uf))(x) $, we have
\begin{align}
    (U^\dagger (\op_h(a) (Uf)))(x) &= \left( \frac{s}{h} \right)^{d/2} \left[  \left( \frac{s}{h} \right)^{d/2} \frac{1}{(2\pi h)^d} \int_{\mathbb{R}^{2d}} e^{i \langle x - \frac{s}{h} z, \xi \rangle / h} a(\xi) f(z) \, dz \, d\xi \right]\\
    &= \left( \frac{s^d}{h^d} \frac{1}{(2\pi h)^d} \right) \iint_{\mathbb{R}^{2d}} e^{i \langle \frac{s}{h} x - \frac{s}{h} z, \xi \rangle / h} a(\xi) f(z) \, dz \, d\xi \\
    &= \frac{1}{\left(2\pi \frac{h^2}{s}\right)^d} \iint_{\mathbb{R}^{2d}} e^{\frac{i}{(h^2/s)} \langle x - y, \xi \rangle } a(\xi) f(z) \, dz \, d\xi \\
    &= \op_{\frac{h^2}{s}} a \quad 
\end{align}
where $z$ is a dummy variable. Since we are given $A_0=\op_h(a)$, we have
    \begin{equation}
    U^\dagger A_0 U = \op_{\frac{h^2}{s}} a.
    \end{equation}

\end{proof}

\REV{
\begin{lemma} \label{lem:UVU}
Suppose that $v(x) \in S(1)$ is a function of only $x \in \mathbb{R}^d$. For $V = \op_h (v)$ and $U$ as defined in \cref{eq:u_def}, we have
    \begin{equation}
        U^\dagger V U 
        = v\left( \frac{s}{h}x \right).
    \end{equation}
\end{lemma}
The proof follows straightforward calculation, which we defer to \cref{sec:appendix} for completeness.
}

\subsection{Proof of \cref{thm:main_S(1)symbol}}\label{sec:subsection_proof_main_thm}

\begin{proof} 
    We introduce the unitary rescaling operator $ U $ defined by
\begin{equation} \label{U_h_0}
U f(x) = \left( \frac{h}{h_0} \right)^{d/2} f \left( \frac{h}{h_0} x \right).
\end{equation}
where we choose $s=h_0$ in Lemma 4. And in this proof, $h_0 = \Delta t$.
\REV{\cref{lem:UVU}} shows that
\begin{equation}
U^\dagger V U f(x) = v \left( \frac{h_0}{h} x \right) f(x),
\end{equation}
and thanks to \cref{lem:UAU}, we also have \begin{equation}
    U^\dagger A_0U=\Op_{\frac{h^2}{s}}a.
\end{equation}
We define $h'=\frac{h^2}{h_0}, \tilde{s} = \frac{s}{h_0}, \tilde{\tau} = \frac{\tau}{h_0}$ and apply $U^\dagger $ to the left and $U$ to the right of the key commutator 
$[ [ V, e^{\frac{i s}{h^2} A_0} V e^{-\frac{i s}{h^2} A_0} ], e^{\frac{i \tau}{h^2} A_0} V e^{-\frac{i \tau}{h^2} A_0} ]$ of interest. Consequently, the following are obtained
\begin{align*}
&U^\dagger [ [ V, e^{\frac{i s}{h^2} A_0} V e^{-\frac{i s}{h^2} A_0} ], e^{\frac{i \tau}{h^2} A_0} V e^{-\frac{i \tau}{h^2} A_0} ] U
\\
=& \left[ [ U^\dagger V U, U^\dagger e^{\frac{i s}{h^2} A_0} U U^\dagger VU U^\dagger e^{-\frac{i s}{h^2} A_0} U ], U^\dagger e^{\frac{i \tau}{h^2} A_0} U U^\dagger VU U^\dagger e^{-\frac{i \tau}{h^2}A_0}U \right] \\
=& \left[ [ \tilde{V}, e^{\frac{i s}{h^2} U^\dagger A_0 U} \tilde{V} e^{-\frac{i s}{h^2} U^\dagger A_0 U}], e^{\frac{i \tau}{h^2} U^\dagger A_0 U} \tilde{V}e^{-\frac{i \tau}{h^2} U^\dagger A_0 U}  \right] \\
=& \left[ [ \tilde{V}, e^{\frac{i}{h'} \tilde{s} \Op_{h'}(a)} \tilde{V} e^{-\frac{i}{h'} \tilde{s} \Op_{h'}(a)}], e^{\frac{i}{h'} \tilde{\tau} \Op_{h'}(a)} \tilde{V} e^{-\frac{i}{h'} \tilde{\tau} \Op_{h'}(a)} \right],
\end{align*}
where $\tilde V: = U^\dagger V U$ and  the last line is obtained using \cref{lem:UAU}.
By \cref{U_h_0} and the fact that $h_0/h = \sqrt{h_0/h'}$, we have \begin{equation}
    \tilde Vf(x)=v\left(\sqrt{\frac{h_0}{h'}}x\right)f(x).
\end{equation}
It is readily observed that $\tilde{V}$ belongs to $\tilde{S}_{\frac{1}{2}}$ with respect to the two parameters $h'$ and $h_0$, both of which range within $(0,1]$, so the Weyl calculus for the symbol class $\tilde{S}_\frac{1}{2}$~\cite{NonnenmacherSjostrandZworski2014} applies. Specifically we apply the Egorov's theorem for $\tilde{S}_\frac{1}{2}$ symbols~\cite[Lemma 3.12]{NonnenmacherSjostrandZworski2014} with the semiclassical parameter $h'$ in our case, we have
\begin{equation}
    e^{\frac{i}{h'} \tilde{s} \Op_{h'}(a)} \tilde{V} e^{-\frac{i}{h'} \tilde{s} \Op_{h'}(a)}=\Op_{h'}(w_{\tilde{s}}), \quad \ e^{\frac{i}{h'} \tilde{\tau} \Op_{h'}(a)} \tilde{V} e^{-\frac{i}{h'} \tilde{\tau} \Op_{h'}(a)}=\Op_{h'}(w_{\tilde{\tau}})
\end{equation}for $w_{\tilde{s}}, w_{\tilde{\tau}}\in\widetilde{S}_{\frac{1}{2}}(1)$ of bounded norm (with the same semiclassical parameters $h'$ and $h_0$). More precisely, one has
$$w_{\tilde{s}} = \kappa_{\tilde{s}} v + hh_0 r, \ w_{\tilde{\tau}} = \kappa_{\tilde{\tau}} v + hh_0 r$$ where $r \in \tilde{S}_\frac{1}{2}$, $\kappa : = \kappa_1$ and $\kappa_t$ is the symplectomorphism generated by the Hamiltonian vector field $H_{a} = \la\dd_{\xi}a,\dd_x\ra-\la\dd_x a,\dd_{\xi}\ra$,
as defined in \cite[Section 3]{NonnenmacherSjostrandZworski2014}. Note that due to the composition rule, each commutator contributes to one order of $h'$ as the underlying Weyl quantization is in terms of the parameter $h'$.

Then by \cite[Lemma 3.2 and Lemma 3.4]{NonnenmacherSjostrandZworski2014} (with the semiclassical parameters where in our case $h'$ correspond to $h$ in the cited work and $h_0$ to $\tilde{h}$), we have 
\begin{equation}  
\left[ [ \tilde{V}, e^{\frac{i}{h'} \tilde{s} \Op_{h'}(a)} \tilde{V} e^{-\frac{i}{h'} \tilde{s} \Op_{h'}(a)}], e^{\frac{i}{h'} \tilde{\tau} \Op_{h'}(a)} \tilde{V} e^{-\frac{i}{h'} \tilde{\tau} \Op_{h'}(a)} \right]
 =i(h'^2)w_\mathrm{comm} + O(h_0h),
\end{equation}
where $w_\mathrm{comm}$ is an $\tilde{S}_{1/2}$ symbol following the composition rules. Following the Calder\'on-Vaillancourt Theorem for $\tilde{S}_{1/2}$ symbols, $w_\mathrm{comm}$ is bounded from $L^2$ to $L^2$ and is $\Or(1)$ (in terms of the semiclassical parameter $h'$).
Note that $h' = h^2/h_0 \leq h_0 = \Delta t$, given the fact that $h\leq h_0 = \Delta t \leq 1$. Thus, we have 
\begin{equation}  
\sup_{\tau, s\in [-\Delta t, \Delta t]}
\left\lVert \left[ \left[ V, e^{i \tau \frac{A_0}{h^2}} V e^{-i \tau \frac{A_0}{h^2}}\right] ,  e^{i s \frac{A_0}{h^2}} V e^{-i s \frac{A_0}{h^2}}  \right] \right\rVert_{L^2(\mathbb{R}^d) \to L^2(\mathbb{R}^d)} \leq C h_0^2 = C \Delta t^2,
\end{equation}
where $C$ is uniform in both $\Delta t$ and $h = 1/N$.
This completes the proof.

\end{proof}

\section{Conclusion and Discussion} \label{sect:conclusion}

In this work, we address the longstanding challenge of extending superconvergence results for unbounded Hamiltonian simulation to discrete spatial discretizations, which is an important aspect of numerical implementations in quantum computing. Using the second-order Magnus expansion, we establish the proof of superconvergence in the discrete setting. This demonstrates that the error constant remains uniform in the number of spatial discretization points $N$. This result validates the robustness of the superconvergence phenomenon for finite $N$.
 
A key innovation of our approach lies in the application of a novel semiclassical analysis framework based on the two-parameter symbol class. We first map the discrete matrix problem to a continuous-variable problem. Notably, a direct application of Egorov's theorem, as used in the continuous case~\cite{AnFangLin2022, FangLiuSarkar2024}, fails here because it requires a time scaling up to $\Or(h^{-2})$, which exceeds the validity range of Egorov's theorem, limited by the Ehrenfest time $\Or(\log(h^{-1}))$.
To overcome this challenge, we apply a unitary transformation that allows us to identify two semiclassical parameters via the spatial grid size and time step size. We then show that the relevant operators, after the transformation, belong to the two-parameter symbol class $\tilde{S}_{1/2}$. The associated symbol calculus guarantees that the error bounds remain uniform in both semiclassical parameters, which in turn translates to uniformity in $N$. To our knowledge, this is the first application of the two-parameter symbol class to numerical analysis, which may be of independent interest.

 \REV{Note that recently (after the appearance of this work), in \cite{FangLiuZhu2025}, one of the authors constructed higher-order quantum Magnus algorithms. These algorithms achieve $\tilde{O}(T\mathrm{polylog}(1/\epsilon))$ cost scaling, comparable to that of the truncated Dyson series. In addition, when truncated to a finite order, the resulting quantum Magnus algorithms are shown to exhibit commutator scaling for general time-dependent Hamiltonians (while truncated Dyson series does not), while still depending on the time derivative of the Hamiltonian only in a polylogarithmic fashion.
Such an extension is nontrivial. The corresponding circuit design and analysis are much more involved: for example, to attain $\mathrm{polylog}(1/\epsilon)$ scaling, both the circuit complexity and analytical bounds must grow at most polynomially in the truncation order $p$, even though the high-order Magnus expansion contains $p!$ terms for a $p$-th order truncation. This challenge is carefully addressed in \cite{FangLiuZhu2025} and a continuous-variable analysis for unbounded case was discussed in \cite{FangZhang2025}. We expect that our discrete superconvergence proof techniques (the two-parameter symbol class and discrete microlocal analysis) can also be extended to these higher-order quantum Magnus circuits, which we leave for future work.}
Extending the superconvergence proof to higher-order Magnus expansions remains an open challenge, particularly given the increased complexity of commutator structures and cancellation mechanisms. Another interesting direction is to explore whether superconvergence behavior can be observed and rigorously established in other systems, such as spin systems, fermionic systems, and bosonic systems. Finally, our semiclassical framework based on two-parameter symbol class may be applied to other numerical discretization methods, including Trotterization and various quantum and classical algorithms. It can serve as a useful and effective tool for establishing algorithmic uniformity concerning discretization parameters.

\bigskip
\section*{Declarations}

\textbf{Funding}. This material is based upon work supported by the U.S. Department of Energy, Office of Science, Accelerated Research in Quantum Computing Centers, Quantum Utility through Advanced Computational Quantum Algorithms, grant no. DE-SC0025572 (D.F.). We also acknowledge the support from NSF DMS-2347791 (J.Z.) and NSF CAREER award DMS-2438074 (D.F.).

\noindent
\textbf{Financial interests}. The authors declare that they have no competing financial interests.

\noindent
\textbf{Data Availability}. This study did not generate or use any datasets. All results presented in this paper are proved without reliance on external data.

\appendix
\section{Auxiliary proofs} \label{sec:appendix}

\begin{proof}[Proof of $U^\dagger U = U U^\dagger$ for \cref{{lem:def_u_unitary}}] \label{pf:appendix_UdagU}
We substitute this into the integral to get
\begin{equation}
\langle Uf, g \rangle_{\REV{L^2}} = \int_{\mathbb{R}^d} \left( \frac{h}{s} \right)^{d/2} f(y) \overline{g\left( \frac{s}{h} y \right)} \left( \frac{s}{h} \right)^d \, dy.
\end{equation}
Simplifying the factors of $ h/s $ and $ s/h $ yields
\begin{equation}
\langle Uf, g \rangle_{\REV{L^2}} = \int_{\mathbb{R}^d} f(y) \left( \frac{s}{h} \right)^{d/2} \overline{g\left( \frac{s}{h} y \right)}  \, dy.
\end{equation}
Thus, by the definition of adjoint $ \langle Uf, g \rangle_{\REV{L^2}} = \langle f, U^\dagger g \rangle_{\REV{L^2}} $, we see that:
\begin{equation}
U^\dagger g(x) = \left( \frac{s}{h} \right)^{d/2} g\left( \frac{s}{h} x \right).
\end{equation}
Now we compute $ U^\dagger U $, applying $ U^\dagger $ to $ Uf $:
\begin{equation}
U^\dagger (Uf)(x) = U^\dagger \left( \left( \frac{h}{s} \right)^{d/2} f\left( \frac{h}{s} x \right) \right).
\end{equation}
Using the formula for $ U^\dagger $, we get
\begin{equation}
U^\dagger (Uf)(x) = \left( \frac{s}{h} \right)^{d/2} \left( \frac{h}{s} \right)^{d/2} f(x) = f(x).
\end{equation}
Thus, $ U^\dagger U = I $. 

Similarly, applying $ U $ to $ U^\dagger g $ give us
\begin{equation}
U(U^\dagger g)(x) = U \left( \left( \frac{s}{h} \right)^{d/2} g\left( \frac{s}{h} x \right) \right).
\end{equation}
We use the definition of $ U $ to obtain
\begin{equation}
U(U^\dagger g)(x) = \left( \frac{h}{s} \right)^{d/2} \left( \frac{s}{h} \right)^{d/2} g(x) = g(x).
\end{equation}
Thus, $ UU^\dagger = I $. 
As a result, we have shown that
\begin{equation}
U^\dagger U = UU^\dagger = I,
\end{equation}
Therefore, $ U $ is a unitary operator.
\end{proof}

\begin{proof}[Proof of \cref{lem:UVU}] 
    From \cref{lem:def_u_unitary}, the operator $U$ as defined in \cref{eq:u_def} is unitary, and its adjoint is therefore given explicitly by
\begin{equation} \label{eq:U_dagger}
    (U^{\dagger}g)(x) = \left(\frac{s}{h} \right)^{d/2} g\left(\frac{s}{h} x \right)
\end{equation}
for all $g\in L^2(\RR^d)$. We use these explicit forms of $U$ and $U^\dagger$ to compute $U^\dagger VU$.

Let $V=\op_h(v)$ with $v=v(x)\in S(1)$ depending only on $x$. Then $V$ acts by pointwise multiplication:
\begin{equation}
    (Vf)(x)=v(x) f(x), \quad \text{for all} \quad f \in L^2(\RR^d).
\end{equation}
We now compute $(U^\dagger VU)f$ for an arbitrary $f\in L^2(\RR^d)$. We first apply $U$ to $f$, and we have
\begin{equation}
     (Uf)(x)=\left(\frac{h}{s} \right)^{d/2}
    f\left(\frac{h}{s} x \right).
\end{equation}
Next, applying $V$ to $Uf$ gives us
\begin{equation}
    (VUf)(x)=v(x)(Uf)(x)=v(x)\left(\frac{h}{s} \right)^{d/2}
    f\left(\frac{h}{s} x \right).
\end{equation}
Finally, we apply $U^\dagger$ to $VUf$, and using the explicit expression of $U^\dagger$ given in \cref{eq:U_dagger}, we obtain:
\begin{equation}\label{eq:UVUf}
     \left(U^\dagger VUf\right)(x)=
    \left(\frac{s}{h} \right)^{d/2}
    \left[v\left(\frac{s}{h} x \right)
        \left(\frac{h}{s} \right)^{d/2}
        f\left(\frac{h}{s} \cdot \frac{s}{h} x \right)\right]\\
\end{equation}
Since $\left(\frac{s}{h} \right)^{d/2}\left(\frac{h}{s} \right)^{d/2}=1$ and $\frac{h}{s}\cdot\frac{s}{h}x = x$, this simplifies to 
\begin{equation} \label{eq:UVUf_simp}
     (U^{\dagger} VU)f(x)= v\left(\frac{s}{h}x \right) f(x).
\end{equation}
Therefore, since \cref{eq:UVUf_simp} holds for all $f\in L^2(\RR^d)$, we have
\begin{equation}
    U^\dagger VU=v\left(\frac{s}{h}x\right).
\end{equation}

\end{proof}

\bibliographystyle{unsrt}
\bibliography{magnus}
\end{document}